\documentclass{amsart}
\usepackage{graphicx}
\usepackage{amssymb}
\usepackage{epstopdf}
\usepackage{nicefrac}
\usepackage{enumerate}
\usepackage{hyperref}
\usepackage{float}
\usepackage{color}
\usepackage{tabularx}

\input xy
\xyoption {all}

\newtheorem{thm}{THEOREM}[section]

\newtheorem{cor}[thm]{COROLLARY}

\newtheorem{defn}[thm]{DEFINITION}
\newtheorem{ex}[thm]{EXAMPLE}

\newtheorem{lemma}[thm]{LEMMA}

\newtheorem{prop}[thm]{PROPOSITION}

\newtheorem{problem}[thm]{PROBLEM}
\newtheorem{remark}[thm]{REMARK}

\newcommand\mor{\mathrel{\stackrel{\makebox[0pt]{\mbox{\normalfont\tiny a}}}{\sim}}}

\begin{document}

\title{Type invariants for  non-abelian odometers}

\author{Steven Hurder}
\address{Steven Hurder, Department of Mathematics, University of Illinois at Chicago, 322 SEO (m/c 249), 851 S. Morgan Street, Chicago, IL 60607-7045}
\email{hurder@uic.edu}

\author{Olga Lukina}
\address{Olga Lukina, Mathematical Institute, 
  Leiden University,
P.O. Box 9512,
2300 RA Leiden,
The Netherlands}
\email{o.lukina@math.leidenuniv.nl}

\thanks{Version date: October 15, 2024}



  \begin{abstract}  
In this work, we introduce the type and typeset   invariants for  equicontinuous group actions on Cantor sets; that is, for generalized odometers. These invariants are collections of equivalence classes of asymptotic Steinitz numbers associated to  the    action.
We   show   the type    is an invariant of the return equivalence class of the action. We introduce the notion of commensurable typesets, and show that two actions which are return equivalent have commensurable typesets.  Examples are given to illustrate the properties of the type and typeset  invariants. The type and typeset invariants are used to define homeomorphism invariants for solenoidal manifolds.
  \end{abstract}

\maketitle

 {   
\section{Introduction}\label{sec-intro}

 The odometer  actions of ${\mathbb Z}$  are  classified  by an infinite sequence of positive integers, up to a suitable notion of equivalence. The  odometer actions of the abelian group ${\mathbb Z}^n$  for $n >1$ cannot be classified by any set of invariants \cite{Thomas2003},  but the works of  Arnold \cite{Arnold1982a,AV1992}, Butler \cite{Butler1965}, Fuchs \cite{Fuchs2015} and Thomas \cite{Thomas2008} for example,  have introduced invariants which can be used for distinguishing ${\mathbb Z}^n$-odometers, and classifying special subclasses of such actions. The   odometer actions of  non-abelian groups are less well  understood, as many interesting new phenomena arise for such actions \cite{HL2019,HL2020,HL2021,HL2022}.
 In this paper, we introduce the \emph{type} and \emph{typeset} invariants for odometers, which are invariants by conjugation and return equivalence. We give examples  that show all possible types can be realized. In the last section of this work, we apply our results to the classification, up to homeomorphism, of solenoidal manifolds. We   first introduce some basic notions.
   
  \subsection{Types} \label{subsec-types}
Let $\vec{m} = \{ m_i \mid 1 \leq i < \infty\}$ be an infinite collection of positive integers. 
 The      \emph{Steinitz number}  (or sometimes called the    \emph{supernatural number}) defined by $\vec{m}$ is  the infinite product
   \begin{equation}\label{eq-steinitzorder}
\xi(\vec{m}) = {\rm lcm} \{ m_1   m_2 \cdots m_{\ell} \mid  \ell > 0\} \ ,
\end{equation}
where ${\rm lcm}$ denotes the least common multiple of the collection of integers.  A Steinitz number $\xi$ can be uniquely written   as the formal product over the set of primes   $\Pi$, 
\begin{equation}
\xi  = \prod_{p \in \Pi} ~ p^{\chi_{\xi}(p)} \  ,
\end{equation}
where the \emph{characteristic function} $\chi_{\xi} \colon \Pi \to \{0,1,\ldots, \infty\}$ counts the multiplicity with which a prime $p$ appears in the infinite product $\xi$.  Associated to a Steinitz number $\xi$ is its  \emph{prime spectrum}.
  \begin{defn}\label{def-primespectrum}
Let  $\Pi = \{2,3,5, \ldots\}$ denote the set of   primes. Given  $\displaystyle \xi = \prod_{p \in \Pi} \ p^{\chi(p)}$, define:
\begin{eqnarray*}
\pi(\xi) ~ & = & ~  \{ p \in \Pi \mid  \chi(p) > 0  \}   \ , ~ \emph{the prime spectrum of} \ \xi  \ ;  \label{eq-primespectrum}\\
\pi_f(\xi) ~ & = & ~  \{ p \in \Pi \mid 0 < \chi(p) < \infty \} \ , ~  \emph{the finite prime spectrum of} \ \xi \ ;  \label{eq-finiteprimespectrum} \\
\pi_{\infty}(\xi) ~ & = & ~  \{ p \in \Pi \mid  \chi(p) = \infty \} \ , ~  \emph{the infinite prime spectrum of} \ \xi \ .\label{eq-infiniteprimespectrum}
\end{eqnarray*}
   \end{defn}

\begin{defn}\label{def-basictype}
Two Steinitz numbers  $\xi$ and $\xi'$ are said to be \emph{asymptotically equivalent} if there exists finite integers $m, m' \geq 1$ such that $m \cdot \xi = m' \cdot \xi'$, and we then write $\xi \mor \xi'$. 
The  asymptotic equivalence class of a Steinitz number $\xi$ is called its  \emph{type}, and denoted by  $\tau[\xi]$.
\end{defn}
   Note that if $\xi \mor \xi'$, then $\pi_{\infty}(\xi) = \pi_{\infty}(\xi')$. The property that $\pi_f(\xi)$ is an \emph{infinite} set is also preserved by asymptotic equivalence of Steinitz numbers, so is an invariant of type.

 The \emph{type} of a Steinitz number was introduced by Baer in \cite[Section 2]{Baer1937}, where the terminology \emph{genus} was used for the type.  Baer used the type    to analyze the classification problem for rank $n$ subgroups  of ${\mathbb Q}^n$,   where the rank of a subgroup  ${\mathcal A} \subset {\mathbb Q}^n$ is the maximum cardinality of a linearly independent subgroup of ${\mathcal A}$.   More detailed discussions can be found   in the works by Arnold \cite[Section 1]{Arnold1982a}),  by   Wilson    \cite[Chapter 2]{Wilson1998},  by  Ribes and   Zalesskii  \cite[Chapter~2.3]{RZ2000} and in Section~\ref{sec-type}.   

\subsection{Cantor actions} \label{subsec-models}
A Cantor action is a group action on a Cantor set by homeomorphisms. Let  $({\mathfrak{X}},\Gamma,\Phi)$   denote a Cantor action  $\Phi \colon \Gamma \times {\mathfrak{X}} \to {\mathfrak{X}}$. We   write $g\cdot x$ for $\Phi(g)(x)$ when appropriate.
The action is \emph{minimal} if  for all $x \in {\mathfrak{X}}$, its   orbit ${\mathcal O}(x) = \{g \cdot x \mid g \in \Gamma\}$ is dense in ${\mathfrak{X}}$.
 The action  $({\mathfrak{X}},\Gamma,\Phi)$ is \emph{equicontinuous} with respect to a metric $d_{\mathfrak{X}}$ on ${\mathfrak{X}}$, if for all ${\varepsilon} >0$ there exists $\delta > 0$, such that for all $x , y \in {\mathfrak{X}}$ and $g \in \Gamma$ we have  that 
 $\displaystyle  d_{\mathfrak{X}}(x,y) < \delta$ implies   $d_{\mathfrak{X}}(g \cdot x, g \cdot y) < {\varepsilon}$.
The property of being equicontinuous    is independent of the choice of the metric   on ${\mathfrak{X}}$ which is   compatible with the topology of ${\mathfrak{X}}$. 
\begin{defn} \cite{CortezPetite2008,CortezMedynets2016} \label{def-odometer}
A   minimal equicontinuous action   $({\mathfrak{X}},\Gamma,\Phi)$ on a Cantor space ${\mathfrak{X}}$ is said to be a $\Gamma$-\emph{odometer}, or just an \emph{odometer} when the group $\Gamma$ is clear from the context.
\end{defn}

A $\Gamma$-odometer has an alternate description in terms of inverse limits of finite actions of $\Gamma$, which is the \emph{algebraic model} for the action. Section~\ref{sec-basics} discusses in further detail the properties of $\Gamma$-odometers.

Let ${\mathcal G} = \{\Gamma = \Gamma_0 \supset \Gamma_1 \supset \Gamma_2 \supset \cdots\}$ be a descending chain of finite index subgroups. Note that the subgroups $\Gamma_{\ell}$ are \emph{not assumed to be normal} in $\Gamma$. 
 Let $X_{\ell} = \Gamma/\Gamma_{\ell}$ and note that  $\Gamma$ acts transitively on the left on the finite set  $X_{\ell}$.    
The inclusion $\Gamma_{\ell +1} \subset \Gamma_{\ell}$ induces a natural $\Gamma$-invariant quotient map $p_{\ell +1} \colon X_{\ell +1} \to X_{\ell}$.
 Introduce the inverse limit 
  \begin{eqnarray} 
X_{\infty} & \equiv &  \lim_{\longleftarrow} ~ \{ p_{\ell +1} \colon X_{\ell +1} \to X_{\ell}  \mid \ell \geq 0 \} \label{eq-invlimspace}\\
& = &  \{(x_0, x_1, \ldots ) \in X_{\infty}  \mid p_{\ell +1 }(x_{\ell + 1}) =  x_{\ell} ~ {\rm for ~ all} ~ \ell \geq 0 ~\} ~ \subset \prod_{\ell \geq 0} ~ X_{\ell} \  .  \nonumber
\end{eqnarray}
Then $X_{\infty}$  is a Cantor space with the Tychonoff topology, where the left actions of $\Gamma$ on the factors $X_{\ell}$ induce    a minimal  isometric action  denoted by  $\Phi_{\infty}  \colon \Gamma \times X_{\infty} \to X_{\infty}$, where $X_{\infty}$ is given the metric induced by the discrete $\Gamma$-invariant metrics on factors $X_{\ell}$. The following is a folklore result:
\begin{thm}\label{thm-gpchain}
Let  $({\mathfrak{X}},\Gamma,\Phi)$ be an odometer, then there exists a group chain  ${\mathcal G}$ and a homeomorphism $\Psi \colon {\mathfrak{X}} \to X_{\infty}$ conjugating the action $\Phi$ with the action $\Phi_{\infty}$.
\end{thm}
The choice of   the odometer $(X_{\infty}, \Gamma, \Phi_{\infty})$ derived from a group chain ${\mathcal G}$ in Theorem~\ref{thm-gpchain} is called an \emph{algebraic model} for $({\mathfrak{X}},\Gamma,\Phi)$. The construction ${\mathcal G}$ uses an adapted neighborhood basis for the odometer, as defined by Definition~\ref{def-adaptednbhds}.
The choice of the group chain  ${\mathcal G}$   is not unique, but two choices are related by a ``tower  equivalence'' of group chains as described in Section~\ref{sec-proofs}. 

 \subsection{Results}\label{subsec-results} 

We associate to a group chain ${\mathcal G}$ a Steinitz number
\begin{equation}\label{eq-orderG}
\xi({\mathcal G}) = {\rm lcm} \{  {\rm Index}[\Gamma : \Gamma_{\ell}]  = \#(\Gamma/\Gamma_{\ell})   \mid  \ell > 0\} \ .
\end{equation}
Our first result is that $\xi({\mathcal G})$ is independent of the choice of the group chain ${\mathcal G}$. In fact, more is true, as the asymptotic equivalence class of $\xi({\mathcal G})$ - its typeset - is invariant under return equivalence. As explained in Section~\ref{sec-basics}, the notion of return equivalence is    the condition that the restrictions of the actions to adapted clopen subsets are isomorphic. It is the group action counterpart of Morita equivalence for groupoids, and is a natural equivalence relation that arises in many contexts.   

  \begin{thm}\label{thm-main2}
  Let    $({\mathfrak{X}},\Gamma,\Phi)$  be an odometer. There is a well-defined Steinitz number  $\xi({\mathfrak{X}},\Gamma,\Phi)$ associated to the action, called its \emph{Steinitz order}.  Moreover, if the $\Gamma'$-odometer  $({\mathfrak{X}}',\Gamma',\Phi')$  is return equivalent to      $({\mathfrak{X}},\Gamma,\Phi)$, then   $\tau[\xi({\mathfrak{X}},\Gamma,\Phi)] = \tau[\xi({\mathfrak{X}'},\Gamma',\Phi')]$; that is, their types are equal.
      \end{thm}
       The type of the action is defined to be  $\tau[\mathfrak{X},\Gamma,\Phi] = \tau[\xi({\mathcal G})]$ for a choice of group chain ${\mathcal G}$.

  The \emph{typeset} invariants of an odometer are a more refined set of invariants that are used, in particular,  to distinguish actions of abelian groups of rank greater than 1. When $\Gamma$ is non-abelian, their definition requires we introduce the normal chain associated to a group chain ${\mathcal G}$.  
  
  Recall that given a subgroup $H \subset G$, the normal core $C(H)  \subset H$ is its largest normal subgroup. 
  If $H$ has finite index, then $C(H)$ also has finite index in $G$.  We associate to $\mathcal G$ the group chain   ${\mathcal G}^c = \{\Gamma = \Gamma_0 \supset C_1 \supset C_2 \supset \cdots\}$ where $C_{\ell} = C(\Gamma_{\ell}) \subset \Gamma_{\ell}$ is its normal core in $\Gamma$.

Let $\gamma \in \Gamma$,  and let $\langle \gamma \rangle \subset \Gamma$ denote the subgroup it generates. For   $\ell > 0$, the intersection 
$\displaystyle \langle \gamma \rangle_{\ell} = \langle \gamma \rangle \cap C_{\ell}$
  is a subgroup of finite index in $\langle \gamma \rangle \cong {\mathbb Z}$. We thus obtain a group chain in  $\langle \gamma \rangle$, denoted
  \begin{equation}\label{eq-gammachain}
{\mathcal C}_{\gamma} = \{ \langle \gamma \rangle \supset \langle \gamma \rangle_1 \supset \langle \gamma \rangle_2 \supset \cdots \} \ .
\end{equation}

\begin{defn}\label{def-typegamma}
Let ${\mathcal G}  = \{\Gamma = \Gamma_0 \supset \Gamma_1 \supset \Gamma_2 \supset \cdots\}$ be a group chain. For $\gamma \in \Gamma$, the type $\tau[\gamma]$ of $\gamma$ is the asymptotic equivalence class of the Steinitz order 
\begin{equation}\label{eq-typegamma}
\xi(\gamma) = {\rm lcm} \{   \#( \langle \gamma \rangle /\langle \gamma \rangle_{\ell}  ) \mid \ell > 0 \}  \ .
\end{equation}
The \emph{typeset} of $\mathcal G$ is the collection
\begin{equation}\label{eq-typeaction}
{\Xi}[{\mathcal G}] = \{\tau[\gamma] \mid \gamma \in \Gamma \} \ .
\end{equation}
\end{defn}

   Note that we allow $\gamma \in \Gamma$ to be the identity above, where $\tau[e] = \{0\}$.

Here are our next two results:

  \begin{thm}\label{thm-main2a}
  Let    $({\mathfrak{X}},\Gamma,\Phi)$  be an odometer, and let $\mathcal G$ be a group chain model. Then the typeset  ${\Xi}[{\mathcal G}]$ is independent of the choice of ${\mathcal G}$.   Thus, ${\Xi}[{\mathcal G}]$ is an invariant under isomorphism of the action.
    \end{thm}  
  The typeset of the action is defined to be  $\Xi[\mathfrak{X},\Gamma,\Phi] = \Xi[{\mathcal G}]$ for a choice of ${\mathcal G}$.
   \begin{thm}\label{thm-main2b}   If the $\Gamma$-odometer  $({\mathfrak{X}},\Gamma,\Phi)$ and   $\Gamma'$-odometer  $({\mathfrak{X}'},\Gamma',\Phi')$ are    return equivalent, then 
  their   typesets  $\Xi[ \mathfrak{X},\Gamma,\Phi]$ and $\Xi[\mathfrak{X}',\Gamma',\Phi']$ are commensurable.
     \end{thm}
     
The notion of commensurable typesets is modeled on the notion of commensurable groups, and given in Definition~\ref{def-commensurable}. The type of an individual $\gamma \in \Gamma$ need not be preserved by return equivalence; rather, it is the collection of all these types - the typeset - that is preserved.

 Various classes of Cantor actions admit stronger results.
For example, we have:
 \begin{cor}\label{cor-abelian}
  Suppose that $\Gamma$ and $\Gamma'$ are   abelian. If  the odometers  $({\mathfrak{X}},\Gamma,\Phi)$  and  $({\mathfrak{X}}',\Gamma',\Phi')$  are return equivalent, then    their   typesets  $\Xi[ \mathfrak{X},\Gamma,\Phi]$ and $\Xi[\mathfrak{X}',\Gamma',\Phi']$ are equal.
 \end{cor}
  Example~\ref{ex-twistedabelian} shows that the typesets need not be equal for return equivalent odometers, if one of the groups is  abelian, and the other is only virtually abelian. Section~\ref{sec-abelian} also gives a selection of ${\mathbb Z}^n$-odometers  to illustrate the definitions of types and typesets. In particular, it is shown that  all types can be realized for all $n\geq 1$.

 A finitely-generated group $\Gamma$ is said to be \emph{renormalizable} if there exists a proper self-embedding $\varphi \colon \Gamma \to \Gamma$ whose image has finite index \cite{HLvL2020}. Another   name for this property is that $\Gamma$ is \emph{finitely non-co-Hopfian}. The embedding $\varphi$ defines a group chain in $\Gamma$ which gives rise to a $\Gamma$-odometer   $(X_{\varphi}, \Gamma, \Phi_{\varphi})$. The properties of the odometers  obtained this way are studied in the work \cite{HLvL2020}. The classification problem for proper  self-embeddings $\varphi \colon {\mathbb Z}^n  \to {\mathbb Z}^n$ was related to the class field theory of number fields  in the   works of Sabitova   \cite{Sabitova2022,Sabitova2024}.

   \begin{thm}\label{thm-renorm}
Let $\varphi \colon \Gamma \to \Gamma$ be a renormalization with associated $\Gamma$-odometer $(X_{\varphi}, \Gamma, \Phi_{\varphi})$. Then  the action has a well-defined type    $\tau[X_{\varphi}, \Gamma, \Phi_{\varphi}]$ and typeset $\Xi [X_{\varphi}, \Gamma, \Phi_{\varphi}]$. Let  $\varphi' \colon \Gamma' \to \Gamma'$ be a renormalization of a second group $\Gamma'$, and assume the odometers $(X_{\varphi}, \Gamma, \Phi_{\varphi})$ and $(X_{\varphi'}, \Gamma', \Phi_{\varphi'})$ are return equivalent, then their types and typesets  are equal.
\end{thm}
 
Many finitely-generated nilpotent groups admit renormalizations, as well as some other classes of  groups, as discussed in \cite{HLvL2020}.   Examples~\ref{ex-trivial} and \ref{ex-2primes}  give   constructions using renormalizable groups which illustrate the conclusions of Theorem~\ref{thm-renorm}. 
For  a finitely-generated, torsion-free nilpotent group $\Gamma$, there are relations between   type  invariants and the dynamics  of  $\Gamma$-odometers  \cite{HL2023}.

Another source of examples arises from group actions on rooted trees, which induce odometer actions on the Cantor sets of ends of the trees.    Section~\ref{subsec-trees} discusses \emph{$d$-regular odometers}, which are odometers that have faithful representations as actions on   $d$-regular trees.

\begin{thm}\label{thm-main4}
A $\Gamma$-odometer $({\mathfrak{X}},\Gamma,\Phi)$ which is  d-regular has finite typeset $\Xi [X_\infty,\Gamma,\Phi]$. 
More precisely, suppose $({\mathfrak{X}},\Gamma,\Phi)$ is isomorphic  to an action on a $d$-ary rooted tree, for some $d \geq 2$. 
Let $P_d$ be the set of distinct prime divisors of the integers   $\{2,\ldots, d\}$, and let $N_d = |P_d|$. Then the cardinality of the typeset satisfies
  \begin{align} \label{eq-upperboundtset0}
  |\Xi [X_\infty,\Gamma,\Phi]| ~ \leq ~ \sum_{k=0}^{N_d} \binom{N_d}{k} = \sum_{k=0}^{N_d} \frac{N_d!}{(N_d - k)! k!} \ .
  \end{align}
  Moreover, each type   $\tau \in \Xi[X_\infty,\Gamma,\Phi]$ is represented by   a Steinitz number $\xi$ with empty finite prime spectrum  $\pi_f(\xi)$, and so $\pi(\xi) = \pi_{\infty}(\xi)$. 
\end{thm}
  Example~\ref{ex-alltypes} gives examples of $d$-regular odometers  which realize the typesets   in Theorem~\ref{thm-main4}.
 
 Finally, we mention two open problems  about the type invariants for odometers.
   \begin{problem}[Realization] \label{prob-realize}
 Given a finitely-generated, torsion-free group $\Gamma$, are there  restrictions  on the types and typesets which can be realized by a  $\Gamma$-odometer?
 \end{problem}
 The work of Arnold \cite{Arnold1982a} implies that every typeset can be realized by a ${\mathbb Z}^n$-odometer. It seems likely that this is also true  when $\Gamma$ is a  nilpotent group. The  solution to Problem~\ref{prob-realize} for a general finitely generated group $\Gamma$   is unknown.
  \begin{problem}[Classification] \label{prob-classify}
 For an odometer $({\mathfrak{X}},\Gamma,\Phi)$ with type   $\tau[{\mathfrak{X}},\Gamma,\Phi]$ and typeset $\Xi[{\mathfrak{X}},\Gamma,\Phi]$, classify   the $\Gamma$-odometers with the same type and typeset.
  \end{problem}
 
 \subsection{Solenoidal manifolds}\label{subsec-solenoids}
The classification of minimal  equicontinuous Cantor actions is closely related to the  classification problem of solenoidal manifolds, and the type invariants defined above yield homeomorphism invariants for solenoidal manifolds.
We briefly recall their definition and state the applications, with details given in Section~\ref{sec-solenoids}.

A \emph{presentation} is a  sequence of \emph{proper finite covering} maps 
${\mathcal P} = \{\, q_{\ell} \colon  M_{\ell} \to M_{\ell -1} \mid  \ell \geq 1\}$, where each $M_{\ell}$ is a compact connected manifold without boundary of dimension $n$.   The inverse limit     
\begin{equation}\label{eq-presentationinvlim}
{\mathcal M}_{{\mathcal P}} \equiv \lim_{\longleftarrow} ~ \{ q_{\ell } \colon M_{\ell } \to M_{\ell -1}\} ~ \subset \prod_{\ell \geq 0} ~ M_{\ell} ~  
\end{equation}
is   the  \emph{weak solenoid}, or \emph{solenoidal manifold},    associated to ${\mathcal P}$. The set ${\mathcal M}_{{\mathcal P}}$ is given  the relative  topology, induced from the product topology, so that ${\mathcal M}_{{\mathcal P}}$ is  compact and connected.   The initial factor $M_0$ is called the \emph{base manifold} of ${\mathcal M}_{{\mathcal P}}$. 

 McCord showed in \cite{McCord1965} that a 
solenoidal manifold  ${\mathcal M}_{{\mathcal P}}$ is a foliated space with foliation ${\mathcal F}_{{\mathcal P}}$, in the sense of \cite{MS2006}, where the leaves of ${\mathcal F}_{{\mathcal P}}$ are coverings of the base manifold $M_0$ via the projection map onto the first factor, ${\widehat{q}}_0 \colon {\mathcal M}_{{\mathcal P}} \to M_0$, restricted to the path-connected components of  ${\mathcal M}_{{\mathcal P}}$. Solenoidal manifolds are     \emph{matchbox manifolds} of dimension $n$ in the terminology of \cite{ClarkHurder2013,CHL2019,HL2019}.  
Solenoidal manifolds have been studied for their geometric properties    \cite{Sullivan2014,Verjovsky2014,Verjovsky2022}, for their analytic and index theory properties  \cite{Connes1994,MS2006}, and  from various number-theoretic viewpoints  \cite{Sabitova2022,Sabitova2024}.
 
Here are two applications of the type invariants. Details and proofs are given in Section~\ref{sec-solenoids}.

 \begin{thm}\label{thm-main11}  
 Associated to a presentation $\mathcal P$ is a well-defined Steinitz number $\xi({\mathcal P})$, and the  type $\tau[\xi({\mathcal P})]$    depends only on the homeomorphism class of ${\mathcal M}_{{\mathcal P}}$. We denote this type by $\tau[{\mathcal M}_{{\mathcal P}}]$.
 \end{thm}

For orientable 1-dimensional solenoids, where each map $q_{\ell} \colon {\mathbb S}^1 \to {\mathbb S}^1$ is orientable,   Bing observed in \cite{Bing1960} that    if  $\tau[{\mathcal M}_{{\mathcal P}}] = \tau[{\mathcal M}_{{\mathcal P}'}]$ then   ${\mathcal M}_{\mathcal P}$ and ${\mathcal M}_{{\mathcal P}'}$  are homeomorphic.  McCord showed in \cite[Section~2]{McCord1965}   the converse, that if  ${\mathcal M}_{\mathcal P}$ and ${\mathcal M}_{{\mathcal P}'}$  are homeomorphic, then  $\tau[{\mathcal M}_{{\mathcal P}}] = \tau[{\mathcal M}_{{\mathcal P}'}]$.   (See also   \cite{AartsFokkink1991}.)  Together these results yield the well-known classification:  
\begin{thm}\label{thm-onedimSol}
For 1-dimensional solenoidal manifolds, 
  ${\mathcal M}_{\mathcal P}$ and ${\mathcal M}_{{\mathcal P}'}$ are homeomorphic  if and only if    $\tau[{\mathcal M}_{\mathcal P}] = \tau[{\mathcal M}_{{\mathcal P}'}]$. 
\end{thm}
For solenoidal manifolds of   dimension $n \geq 2$, no such classification by invariants can exist, even for the case when  the base manifold $M_0 = {\mathbb T}^n$ (see Thomas \cite{Thomas2003}.) The most one can hope for is to define invariants which distinguish various homeomorphism classes.
Associated to a presentation $\mathcal P$ is a well-defined typeset $\Xi({\mathcal P})$, which provides another       invariant of the homeomorphism class of ${\mathcal M}_{\mathcal P}$.

 \begin{thm}\label{thm-main21}  
 Assume that ${\mathcal M}_{{\mathcal P}}$ and ${\mathcal M}_{{\mathcal P}'}$ are homeomorphic solenoidal manifolds, then the typesets $\Xi[ {{\mathcal P}}]$ and $\Xi[ {{\mathcal P}'}]$  are commensurable.
 \end{thm}
 
 \subsection{Structure}
 The remainder of this paper proceeds as following.
 
   Section~\ref{sec-basics}  discusses the basic results about odometers required.

   Section~\ref{sec-type}   discusses the basic results about types and typesets required.

   Section~\ref{sec-proofs}   gives the proofs of the results above about classifying odometers.

 Section~\ref{sec-abelian}  gives some basic examples of odometers and calculates their types and typesets.

 Section~\ref{sec-regular}  discusses $d$-regular odometers and calculates their types and typesets.

Section~\ref{sec-solenoids}  discusses solenoidal manifolds and their types and typesets.

 \section{Cantor actions}\label{sec-basics}

We recall some of the basic 
 properties of     Cantor actions, as required for the proofs of the results in Section~\ref{sec-intro}.
More complete discussions of the properties of equicontinuous Cantor actions are given in     the text by Auslander \cite{Auslander1988}, the papers by Cortez and Petite  \cite{CortezPetite2008}, Cortez and Medynets  \cite{CortezMedynets2016},    and   the authors' works, in particular   \cite{DHL2016a,DHL2016c} and   \cite[Section~3]{HL2021}.

\subsection{Basic concepts}\label{subsec-basics}

Assume that ${\mathfrak{X}}$ is a Cantor space. 
Let ${\rm CO}({\mathfrak{X}})$ denote the collection  of all clopen (closed and open) subsets of  ${\mathfrak{X}}$, which forms a basis for the topology of ${\mathfrak{X}}$. 
For $\phi \in {\rm Homeo}({\mathfrak{X}})$ and    $U \in {\rm CO}({\mathfrak{X}})$, the image $\phi(U) \in {\rm CO}({\mathfrak{X}})$.  Recall that the action $({\mathfrak{X}},\Gamma,\Phi)$ is \emph{minimal} if  for all $x \in {\mathfrak{X}}$, its   orbit ${\mathcal O}(x) = \{g \cdot x \mid g \in \Gamma\}$ is dense in ${\mathfrak{X}}$.

The following   result is folklore, and a proof is given in \cite[Proposition~3.1]{HL2020}.
 \begin{prop}\label{prop-CO}
 For ${\mathfrak{X}}$ a Cantor space, a minimal   action   $\Phi \colon \Gamma \times {\mathfrak{X}} \to {\mathfrak{X}}$  is  equicontinuous  if and only if  the $\Gamma$-orbit of every $U \in {\rm CO}({\mathfrak{X}})$ is finite for the induced action $\Phi_* \colon \Gamma \times {\rm CO}({\mathfrak{X}}) \to {\rm CO}({\mathfrak{X}})$.
\end{prop}
 That is, given a clopen set $U \in {\rm CO}({\mathfrak{X}})$,  $\Gamma$ acts on the orbits of $U$ by finite permutations. This justifies saying that a minimal equicontinuous Cantor action of $\Gamma$ is a (generalized) odometer.
   
We say that $U \subset {\mathfrak{X}}$  is \emph{adapted} to the action $({\mathfrak{X}},\Gamma,\Phi)$ if $U$ is a   \emph{non-empty clopen} subset, and for any $g \in \Gamma$, 
if $\Phi(g)(U) \cap U \ne \emptyset$ implies that  $\Phi(g)(U) = U$.   Given  $x \in {\mathfrak{X}}$ and clopen set $x \in W$, there is an adapted clopen set $U$ with $x \in U \subset W$. (For example, see the proof of   \cite[Proposition~3.1]{HL2020}.)

For an adapted set $U$,   the set of ``return times'' to $U$, 
 \begin{equation}\label{eq-adapted}
\Gamma_U = \left\{g \in \Gamma \mid g \cdot U  \cap U \ne \emptyset  \right\}  
\end{equation}
is a subgroup of   $\Gamma$, called the \emph{stabilizer} of $U$.      
  Then for $g, g' \in \Gamma$ with $g \cdot U \cap g' \cdot U \ne \emptyset$ we have $g^{-1} \, g' \cdot U = U$, hence $g^{-1} \, g' \in \Gamma_U$. Thus,  the  translates $\{ g \cdot U \mid g \in \Gamma\}$ form a finite clopen partition of ${\mathfrak{X}}$, and are in 1-1 correspondence with the quotient space $X_U = \Gamma/\Gamma_U$. Then $\Gamma$ acts by permutations of the finite set $X_U$ and so the stabilizer group $\Gamma_U \subset G$ has finite index.  Note that this implies that if $V \subset U$ is a proper inclusion of adapted sets, then the inclusion $\Gamma_V \subset \Gamma_U$ is also proper.

\begin{defn}\label{def-adaptednbhds}
Let  $({\mathfrak{X}},\Gamma,\Phi)$   be an odometer.
A properly descending chain of clopen sets, ${\mathcal U} = \{U_{\ell} \subset {\mathfrak{X}}  \mid \ell > 0\}$, is said to be an \emph{adapted neighborhood basis} at $x \in {\mathfrak{X}}$ for the action $\Phi$,   if
    $x \in U_{\ell +1} \subset U_{\ell}$  is a proper inclusion for all $ \ell > 0$, with     $\cap_{\ell > 0}  \ U_{\ell} = \{x\}$, and  each $U_{\ell}$ is adapted to the action $\Phi$.
\end{defn}
Given $x \in {\mathfrak{X}}$ and   ${\varepsilon} > 0$, Proposition~\ref{prop-CO} implies there exists an adapted clopen set $U \in {\rm CO}({\mathfrak{X}})$ with $x \in U$ and ${\rm diam}(U) < {\varepsilon}$.  Thus, one can choose a descending chain ${\mathcal U}$ of adapted sets in ${\rm CO}({\mathfrak{X}})$ whose intersection is $x$, from which the following result follows:

\begin{prop}\label{prop-adpatedchain}
Let  $({\mathfrak{X}},\Gamma,\Phi)$   be an odometer. Given $x \in {\mathfrak{X}}$, there exists an adapted neighborhood basis ${\mathcal U}$ at $x$ for the action $\Phi$.
 \end{prop}
\begin{cor}\label{cor-Uchain}  
Let  $({\mathfrak{X}},\Gamma,\Phi)$   be an odometer, and   ${\mathcal U}$ be an adapted neighborhood basis. Set  $\Gamma_{\ell} = \Gamma_{U_{\ell}}$,  with $\Gamma_0 = \Gamma$, then  $\displaystyle  {\mathcal G}_{{\mathcal U}} = \{\Gamma_0 \supset \Gamma_1 \supset \cdots \}$ is a  descending chain of finite index subgroups.
 \end{cor}

\subsection{Equivalence of Cantor actions}\label{subsec-equivalence}

We next recall the   notions of equivalence of  Cantor actions  used in this work. 
The first and strongest   is that  of 
  {isomorphism}, which   is a   generalization  of the usual notion of conjugacy of topological actions.
  The definition below agrees with the usage in the papers     \cite{CortezMedynets2016,HL2020,Li2018}.

 \begin{defn} \label{def-isomorphism}
Cantor actions $({\mathfrak{X}}_1, \Gamma_1, \Phi_1)$ and $({\mathfrak{X}}_2, \Gamma_2, \Phi_2)$ are said to be \emph{isomorphic}  if there is a homeomorphism $h \colon {\mathfrak{X}}_1 \to {\mathfrak{X}}_2$ and group isomorphism $\Theta \colon \Gamma_1 \to \Gamma_2$ so that 
\begin{equation}\label{eq-isomorphism}
\Phi_1(g) = h^{-1} \circ \Phi_2(\Theta(g)) \circ h   \in   {\rm Homeo}({\mathfrak{X}}_1) \   \textrm{for  all} \ g \in \Gamma_1 \ .
\end{equation}
 \end{defn}

The notion of \emph{return equivalence} for odometers is  weaker than the notion of isomorphism, and is natural when considering the odometers defined by the holonomy actions for solenoidal manifolds, as considered in   the works  \cite{HL2019,HL2020,HL2021}, and later in Section~\ref{sec-proofs}.

For an odometer $({\mathfrak{X}}, \Gamma, \Phi)$ and   an adapted set $U \subset {\mathfrak{X}}$, by an   abuse of  notation, we use $\Phi_U$ to denote both the restricted action $\Phi_U \colon \Gamma_U \times U \to U$, and the induced quotient action $\Phi_U \colon {\mathcal H}_U \times U \to U$, where  ${\mathcal H}_U = \Phi(\Gamma_U) \subset {\rm Homeo}(U)$. 
 Then $(U, {\mathcal H}_U, \Phi_U)$ is called the \emph{restricted holonomy action} for $\Phi$, in analogy with the case where $U$ is a transversal to a solenoidal manifold, and 
  ${\mathcal H}_U$ is the holonomy group for this transversal. A technical issue that often arises though, is that while $\Gamma_U \subset \Gamma$ has finite index, the action map $\Phi_U \colon \Gamma_U \to {\mathcal H}_U$ need not be injective, and can in fact can have a large kernel as is the case for example for the actions of weakly branch groups (see Example~\ref{ex-branch}).

   \begin{defn}\label{def-return}
Odometers $({\mathfrak{X}}, \Gamma, \Phi)$ and $({\mathfrak{X}}', \Gamma', \Phi')$ are  \emph{return equivalent} if there exists 
  an adapted set $U \subset {\mathfrak{X}}$ for the action $\Phi$,   and  
  an adapted set $U' \subset {\mathfrak{X}}'$ for the action $\Phi'$,
such that   the  restricted actions $(U, {\mathcal H}_{U}, \Phi_{U})$ and $(U', {\mathcal H}'_{U'}, \Phi'_{U'})$ are isomorphic.
\end{defn}
 Note that if we take $U = {\mathfrak{X}}$ and $U' = {\mathfrak{X}}'$ in Definition~\ref{def-return}, then return equivalence may still be weaker than isomorphism in  Definition~\ref{def-isomorphism}, unless the actions $\Phi$ and $\Phi'$ are topologically free \cite{HL2020,Li2018,Renault2008}.
 
 \subsection{Algebraic Cantor  actions}\label{subsec-gchains}
 
    There is a natural basepoint $x_{\infty} \in X_{\infty}$ given by the cosets of the identity element $e \in \Gamma$, so $x_{\infty} = (e \Gamma_{\ell})$. An adapted neighborhood basis   of $x_{\infty}$ is given by the clopen sets 
\begin{equation}\label{eq-openbasis}
U_{\ell} = \left\{ x = (x_{i}) \in X_{\infty}   \mid  x_i = e \Gamma_i \in X_i~, ~ 0 \leq i \leq \ell ~  \right\} \subset X_{\infty} \ .
\end{equation}
Then there is the tautological identity $\Gamma_{\ell} = \Gamma_{U_{\ell}}$.

Suppose that we are given an odometer $({\mathfrak{X}},\Gamma,\Phi)$, and   an adapted neighborhood basis ${\mathcal U}$. Define subgroups $\Gamma_{\ell} = \Gamma_{U_{\ell}}$,  with $\Gamma_0 = \Gamma$, which form the group chain 
$\displaystyle  {\mathcal G}_{{\mathcal U}} = \{\Gamma_0 \supset \Gamma_1 \supset \cdots \}$. Then we have the folklore result:
\begin{thm}\label{thm-algmodel}
Let $({\mathfrak{X}},\Gamma,\Phi)$ be an odometer, and ${\mathcal U}$ an adapted neighborhood basis. The action $({\mathfrak{X}},\Gamma,\Phi)$ is isomorphic  to the   odometer $(X_{\infty}, \Gamma, \Phi_{\infty})$ constructed from the group chain ${\mathcal G}_{{\mathcal U}}$.
\end{thm}

  \begin{cor}\label{cor-equivalent}
  Let $({\mathfrak{X}}, \Gamma, \Phi)$ be an odometer, and assume that  ${\mathcal G}_{{\mathcal U}}$ and  ${\mathcal G}'_{{\mathcal U}'}$ are adapted neighborhood bases for the action. Then the corresponding algebraic models of the action, $(X_{\infty}, \Gamma, \Phi_{\infty})$ and $(X'_{\infty}, \Gamma, \Phi_{\infty}')$, are isomorphic in the sense of Definition~\ref{def-isomorphism} with $\Theta \colon \Gamma \to \Gamma$ the identity map.
  \end{cor}

\section{Type and typeset} \label{sec-type} 

 The notion of type was introduced in 1937 by   Baer in \cite[Section 2]{Baer1937} as part of the study of the classification problem for   rank $n$  subgroups of ${\mathbb Q}^n$. The work of   Butler \cite{Butler1965} introduced a restricted class of subgroups in ${\mathbb Q}^n$ now called \emph{Butler groups}. The  classification theory for Butler groups was further developed by    Richman \cite{Richman1983} and    Mutzbauer \cite{Mutzbauer1983}, and the   works of   Arnold  (see  \cite[Section 1]{Arnold1982a}),  and Arnold and Vinsonhaler \cite{AV1992}. For a comprehensive treatment of these ideas, see the monograph by  Fuchs \cite{Fuchs2015}.    Thomas applied the type invariants    in his analysis of the classification complexity of these groups in the work   \cite[Section~3]{Thomas2008}.   The applications of type  invariants to profinite groups are discussed in the works by   Ribes \cite[Chapter~1,  Section~4]{Ribes1970}, Wilson    \cite[Chapter 2]{Wilson1998} and   Ribes and   Zalesskii  \cite[Chapter~2.3]{RZ2000}. 
 In this section, we recall     basic notions and properties of Steinitz numbers and their types,  and their definitions     for $\Gamma$-odometers. 
 
\subsection{Types and typesets}\label{subsec-types2}
Recall that a  Steinitz number $\xi$ can be written uniquely as the formal product over the set of primes, 
\begin{equation}\label{eq-defsteinitz}
\xi  = \prod_{p \in \Pi} ~ p^{\chi_{\xi}(p)} \  ,
\end{equation}
where the \emph{characteristic function} $\chi_{\xi} \colon \Pi \to \{0,1,\ldots, \infty\}$ counts the multiplicity with which a prime $p$ appears in the infinite product $\xi$. 
Note that multiplication of Steinitz numbers corresponds to  the sum of their characteristic functions. 

Recall from Definition \ref{def-basictype} that two Steinitz numbers  $\xi$ and $\xi'$ are said to be \emph{asymptotically equivalent} if there exists finite integers $m, m' \geq 1$ such that $m \cdot \xi = m' \cdot \xi'$, and we then write $\xi \mor \xi'$. Recall that the type $\tau[\xi]$ associated to a Steinitz number  $\xi$  is the asymptotic equivalence class   of $\xi$.

 \begin{lemma}\label{lem-charequic}
 $\xi$ and $\xi'$ satisfy $\xi \mor \xi'$ if and only if their characteristic functions $\chi_1, \chi_2$ satisfy
 \begin{itemize}
\item $\chi_1(p) = \chi_2(p)$ for all but finitely many primes $p \in \Pi$, and
\item $\chi_1(p) = \infty$ if and only if $\chi_2(p) = \infty$ for all primes $p \in \Pi$.
\end{itemize}
 \end{lemma}
Given two types  $\tau$ and $\tau'$, we write $\tau \leq \tau'$  if there exists representatives $\xi \in \tau$ and $\xi' \in \tau'$ such that their characteristic functions satisfy $\chi_{\xi}(p) \leq \chi_{\xi'}(p)$ for all primes $p \in \Pi$. Then  two Steinitz numbers $\xi$ and $\xi'$ are asymptotically equivalent if and only if $\chi_{\xi} \leq \chi_{\xi'}$ and $\chi_{\xi'} \leq \chi_{\xi}$. 
\begin{defn}\label{def-typeset}
 A \emph{typeset} $\Xi$ is   a collection of types. 
\end{defn}

There are three operations on types $\tau$ and $\tau'$: product, join and intersection. Let $\chi$ (respectively $\chi'$) be the characteristic function for a representative  $\xi \in \tau$ (respectively $\xi' \in \tau'$),  then:

$\begin{array}{lccccc}
{\rm Product:} &  \tau \cdot \tau' & {\rm type ~ defined ~ by} & \chi^{+}(p)  & = &   \chi(p) + \chi'(p)    \\
{\rm Join:} &  \tau \vee \tau' & {\rm type ~ defined ~ by} & \chi^{\vee}(p)  & = & \max \{ \chi(p) , \chi'(p)\}    \\
{\rm Intersection:} &   \tau \wedge \tau'   & {\rm type ~ defined ~  by} & \chi^{\wedge}(p)  & = &\min \{ \chi(p) , \chi'(p)\}  
\end{array}
$

Note that    $\tau \wedge \tau' \leq \tau \vee \tau' \leq   \tau \cdot \tau'$.
 A  typeset $\Xi$  need not be closed under   the   operations of product, join or intersection. However,    a typeset  $\Xi$ always admits a partial ordering.

 \subsection{Type  for  profinite groups}

A topological group ${\mathfrak{G}}$ is said to be profinite if it is isomorphic to the  inverse limit of finite groups (see    \cite[Chapter 2]{Wilson1998} or \cite[Chapter~2.3]{RZ2000}).

 The \emph{Steinitz order} $\Pi[{\mathfrak{G}}]$ of a   profinite group ${\mathfrak{G}}$ is  defined by the supernatural number associated to a presentation of ${\mathfrak{G}}$, defined as follows. 
For a profinite group ${\mathfrak{G}}$, an open subgroup ${\mathfrak{U}} \subset {\mathfrak{G}}$ has finite index  \cite[Lemma 2.1.2]{RZ2000}. Let ${\mathfrak{D}} \subset {\mathfrak{G}}$ be a closed subgroup, and ${\mathfrak{N}} \subset {\mathfrak{G}}$ is an open normal subgroup, then ${\mathfrak{N}} \cdot {\mathfrak{D}}$ is   an open subgroup of ${\mathfrak{G}}$ and ${\mathfrak{N}} \cap {\mathfrak{D}}$ is an open normal subgroup of ${\mathfrak{D}}$. 

\begin{defn}\label{def-steinitzorderaction}
Let ${\mathfrak{D}} \subset {\mathfrak{G}}$ be a closed subgroup of the profinite group ${\mathfrak{G}}$. Define   Steinitz orders    as follows:
 \begin{enumerate}
\item $\xi({\mathfrak{G}}) =  {\rm lcm} \{\# \ {\mathfrak{G}}/{\mathfrak{N}}  \mid {\mathfrak{N}} \subset {\mathfrak{G}}~ \text{open normal subgroup}\}$,
\item $\xi({\mathfrak{D}}) =  {\rm lcm} \{\# \ {\mathfrak{D}}/({\mathfrak{N}} \cap {\mathfrak{D}}) \mid {\mathfrak{N}} \subset {\mathfrak{G}}~ \text{open normal subgroup}\}$, 
\item $\xi({\mathfrak{G}} : {\mathfrak{D}}) =  {\rm lcm} \{\# \ {\mathfrak{G}}/({\mathfrak{N}} \cdot {\mathfrak{D}})  \mid  {\mathfrak{N}} \subset {\mathfrak{G}}~ \text{open normal subgroup}\}$.
\end{enumerate}
\end{defn}
The Steinitz number $\xi({\mathfrak{G}} : {\mathfrak{D}})$ is called the \emph{relative Steinitz order} of the pair $({\mathfrak{G}}, {\mathfrak{D}} )$.

The Steinitz orders  satisfy the Lagrange identity, where the multiplication is taken in the sense of supernatural numbers  (see \cite{RZ2000,Wilson1998}), and we have
 \begin{equation}\label{eq-productorders}
\xi({\mathfrak{G}})  = \xi({\mathfrak{G}} : {\mathfrak{D}}) \cdot \xi({\mathfrak{D}})  \ .
\end{equation}
 
 In particular,  we always have $\xi[{\mathfrak{D}}] \leq \xi[{\mathfrak{G}}]$. 

 \begin{defn}\label{def-reltype}
  The \emph{relative type} $\Xi[{\mathfrak{G}} : {\mathfrak{D}}]$ is the type of   $\xi({\mathfrak{G}} : {\mathfrak{D}})$, where   ${\mathfrak{D}} \subset {\mathfrak{G}}$ is a closed subgroup of the profinite group ${\mathfrak{G}}$.

 \end{defn}

 \subsection{Type  for  odometers}
 
 Let $(X_{\infty}, \Gamma, \Phi_{\infty})$ be an odometer, defined by a group chain ${\mathcal G} = \{\Gamma = \Gamma_0 \supset \Gamma_1  \supset \cdots\}$.
 Recall that the  \emph{normal core} of $\Gamma_{\ell}$  is the largest normal subgroup   $C_{\ell} \subset \Gamma_{\ell}$. 
The    profinite group $\widehat{\Gamma}({\mathcal G})$ associated to the action $(X_{\infty}, \Gamma, \Phi_{\infty})$ defined by ${\mathcal G}$ is given by:
\begin{defn}\label{def-profinitecomp}
Let ${\mathcal G}$ be a group chain with associated normal core chain $\{C_{\ell} \mid \ell \geq 0\}$, then
\begin{equation}\label{eq-profinitelim}
\widehat{\Gamma}({\mathcal G}) = \lim_{\longleftarrow} \{ \Gamma/C_{\ell +1} \to \Gamma/C_{\ell} \mid \ell \geq 0 \} \ .
\end{equation}
\end{defn}
 
The profinite group $\widehat{\Gamma}({\mathcal G})$ is   a quotient of the full profinite completion of $\Gamma$, but it is typically not equal to the profinite completion. The properties of the group $\widehat{\Gamma}({\mathcal G})$ associated to a Cantor action are studied extensively in the authors' works \cite{DHL2016a,DHL2016c,HL2019,HL2021}.
 The Steinitz order   of $\widehat{\Gamma}({\mathcal G})$ is well-defined, and given by 
 \begin{equation}
\xi(\widehat{\Gamma}({\mathcal G})) =   {\rm lcm} \{   \# ( \Gamma/C_{\ell}) \mid \ell > 0 \}  \ .
\end{equation}
 
  The   type $\tau[\widehat{\Gamma}({\mathcal G})] =   \tau[\xi(\widehat{\Gamma}({\mathcal G}))]$   need not be an invariant of return equivalence of the odometer $(X_{\infty}, \Gamma, \Phi_{\infty})$, as explained below.  

\begin{defn}\label{def-typegroup}
Let $(X_{\infty}, \Gamma, \Phi_{\infty})$ be an odometer.   The type $\tau[X_{\infty}, \Gamma, \Phi_{\infty}]$ of the action is the equivalence class of the Steinitz number
\begin{equation}\label{eq-actionsteinitzorder}
\xi(X_{\infty}, \Gamma, \Phi_{\infty}) = {\rm lcm} \{ \# X_{\ell} = \#(\Gamma/\Gamma_{\ell}) \mid \ell > 0 \}  \ .
\end{equation}
\end{defn}

 Recall that there is a    transitive action     $\widehat{\Phi}_{\infty} \colon \widehat{\Gamma}({\mathcal G}) \times X_{\infty} \to X_{\infty}$ induced by odometer $(X_{\infty}, \Gamma, \Phi_{\infty})$.
The  action $(X_{\infty}, \widehat{\Gamma}({\mathcal G}), \widehat{\Phi}_{\infty})$ is   free precisely 
when the isotropy subgroup ${\mathcal D}({\mathcal G}) \subset  \widehat{\Gamma}({\mathcal G})$ of the action   at the basepoint   $x_{\infty} \in X_{\infty}$ is trivial. In this case, we have a homeomorphism $X_{\infty} \cong \widehat{\Gamma}({\mathcal G})$ that commutes with the action, and so $X_{\infty}$ inherits the structure of a Cantor group from the action.

However, when the action $(X_{\infty}, \widehat{\Gamma}({\mathcal G}), \widehat{\Phi}_{\infty})$ is not free, then ${\mathfrak{D}}({\mathcal G})$ is not trivial, and  we have:
\begin{prop}\label{prop-relsteinitz}
 Let $(X_{\infty}, \Gamma, \Phi_{\infty})$ be an odometer. Then  $\xi(X_{\infty}, \Gamma, \Phi_{\infty}) = \xi(\widehat{\Gamma}({\mathcal G}) : {\mathfrak{D}}({\mathcal G}))$.
 \end{prop}
We omit the proof of this, as it  is a direct consequence of the definitions, and the result is not needed for the proofs of our main theorems.  However, Proposition~\ref{prop-relsteinitz} provides some insights on the properties of the type $\tau[X_{\infty}, \Gamma, \Phi_{\infty}]$ of an odometer with respect to return equivalence.

  The   Lagrange Theorem for profinite groups \eqref{eq-productorders}   implies  that $\tau[\widehat{\Gamma}({\mathcal G})] = \tau[X_{\infty}, \Gamma, \Phi_{\infty}] \cdot  \tau[{\mathcal D}({\mathcal G})]$. The type  $\tau[{\mathcal D}({\mathcal G})]$ need not be invariant under restriction to adapted subsets, and so   $\tau[{\mathcal D}({\mathcal G})]$ need not be invariant under the relation of return equivalence, and thus the same is true for $\tau[\widehat{\Gamma}({\mathcal G})]$. On the other hand,  the proof of Theorem~\ref{thm-main2} shows that  the relative type $\tau[\widehat{\Gamma}({\mathcal G}) : {\mathcal D}({\mathcal G})] = \tau[{\mathcal G}]$   is invariant under restriction to adapted subsets.

  \subsection{Typesets for  odometers}\label{subsec-commensurable}

We next consider the properties of typeset  under restriction, which leads to the   notion of commensurable typesets. Let ${\mathcal G}  = \{\Gamma = \Gamma_0 \supset \Gamma_1 \supset \Gamma_2 \supset \cdots\}$ be a group chain, and  $(X_{\infty}, \Gamma, \Phi_{\infty})$ the associated odometer.

 Let   $H \subset \Gamma$ be a subgroup of finite index such that $\Gamma_{\ell} \subset H$ for some $\ell > 0$.  By omitting some initial terms, we can assume that $\Gamma_1 \subset H$. Then for  each $\ell > 0$, define
\begin{equation}\label{eq-Hcore}
C_{\ell}^H = \bigcap_{\delta \in H} ~  \delta  \Gamma_{\ell} \delta^{-1}    \ .
\end{equation}
Then $C_{\ell}^H \subset \Gamma_{\ell}$ is a subgroup of finite index, as $\Gamma_{\ell}$ has finite index in $\Gamma$. In particular, 
 $C_{\ell} = C_{\ell}^{\Gamma}$ is the normal core of $\Gamma_{\ell}$, and there is an inclusion  $C_{\ell} \subset C_{\ell}^H$ for any choice of $H$.  
  Also note  that   $C_{\ell+1}^H \subset C_{\ell}^H$ for all $\ell \geq 0$.
  
   Given $\gamma \in \Gamma$ define $C_{\gamma, \ell}^H = \langle \gamma \rangle \cap C_{\ell}^H$, 
 and let 
 ${\mathcal C}_{\gamma}^H = \{\langle \gamma \rangle \cap H = C_{\gamma, 0}^H \supset  C_{\gamma, 1}^H  \supset  C_{\gamma, 2}^H \supset \cdots \}$ denote the resulting subgroup chain in $\langle \gamma \rangle$.

\begin{defn}\label{def-Htypegamma}
Given a group chain ${\mathcal G}$ and   $H \subset \Gamma$ a subgroup with $\Gamma_{\ell} \subset H$ for some $\ell \geq 1$, the \emph{$H$-restricted order} of $\gamma \in \Gamma$ is the  Steinitz order with respect to the chain ${\mathcal C}_{\gamma}^H$, 
\begin{equation}\label{eq-Htypegamma}
\xi^H(\gamma) = {\rm lcm} \{   \#( \langle \gamma \rangle /C_{\gamma, \ell}^H  ) \mid \ell > 0 \}  \ .
\end{equation}
The $H-$\emph{restricted typeset} for the chain ${\mathcal G}$  is the collection
\begin{equation}\label{eq-Htypeaction}
\Xi_H[{\mathcal G}] = \{ \tau[\xi^H(\gamma)] \mid \gamma \in \Gamma\} \ .
\end{equation}
\end{defn}
Note that when $H = \Gamma$ we recover  the typeset $\Xi[{\mathcal G}]$ in Definition~\ref{def-typegamma}.

 We can now formulate the   notion of commensurable typesets for group chains ${\mathcal G}$ and ${\mathcal G}'$.  
\begin{defn}\label{def-commensurable}
Let ${\mathcal G} = \{\Gamma = \Gamma_0 \supset \Gamma_1 \supset \Gamma_2 \supset \cdots\}$ and ${\mathcal G}' = \{\Gamma' = \Gamma'_0 \supset \Gamma'_1 \supset \Gamma'_2 \supset \cdots\}$ be group chains. We say that their typesets $\Xi[{\mathcal G}]$ and $\Xi[{\mathcal G}']$ are \emph{commensurable} if there exists a finite index subgroup $ H \subset \Gamma$ with $\Gamma_{\ell} \subset H$ for some $\ell \geq 1$, and $H' \subset \Gamma'$ with $\Gamma'_{\ell'} \subset H'$ for some $\ell' \geq 1$, such that $\Xi_H[{\mathcal G}] = \Xi_{H'}[{\mathcal G}']$.
\end{defn}   
 
 Finally, we give an elementary result about the types of elements $\gamma \in \Gamma$ that follows from basic group theory.
 For all $\ell > 0$, we have that  $\langle \gamma \rangle _{\ell}   \subset  C_{\ell}$,  
 so  the order of the subgroup $\langle \gamma \rangle/\langle \gamma \rangle_{\ell}$ divides the order of $\Gamma/C_{\ell}$  by Lagrange's Theorem. Thus we have: 
\begin{prop}
For a group chain ${\mathcal G}$,  $\tau[{\mathcal G}] \leq \tau[\widehat{\Gamma}({\mathcal G})]$, and for each $\gamma \in \Gamma$, $\tau[\gamma] \leq \tau[\widehat{\Gamma}({\mathcal G})]$.
\end{prop}

\section{Invariance of types} \label{sec-proofs} 

In this section, we give the proofs of the results in Section~\ref{sec-intro} on odometers. The discussion of solenoidal manifolds and the proofs of Theorems~\ref{thm-main11} and \ref{thm-main21} are given in Section~\ref{sec-solenoids}.

The following  notation is used in the following.
Let  $({\mathfrak{X}}, \Gamma, \Phi)$ and $({\mathfrak{X}'}, \Gamma', \Phi')$ denote odometers, and let     ${\mathcal U}$ for $({\mathfrak{X}}, \Gamma, \Phi)$ be a choice of an adapted neighborhood basis for  $({\mathfrak{X}}, \Gamma, \Phi)$, and 
  ${\mathcal U}'$   a choice of an adapted neighborhood basis for $({\mathfrak{X}'}, \Gamma', \Phi')$. In some cases, we have ${\mathfrak{X}} = {\mathfrak{X}'}$.

Let  ${\mathcal G}_{{\mathcal U}}$ be the group chain associated to $\mathcal U$ as in Corollary~\ref{cor-Uchain}, and let  $(X_{\infty}, \Gamma, \Phi_{\infty})$ denote the algebraic model for $({\mathfrak{X}}, \Gamma, \Phi)$ that it determines as in Section~\ref{subsec-models}.
Similarly, let ${\mathcal G}'_{{\mathcal U}'}$  be the group chain associated to ${\mathcal U}'$ and  $(X_{\infty}', \Gamma', \Phi_{\infty}')$ the algebraic model for  $({\mathfrak{X}'}, \Gamma', \Phi')$.

\subsection{Proof of Theorem~\ref{thm-main2}}\label{subsec-proof-main2}

We must show that the Steinitz number  $\xi({\mathfrak{X}},\Gamma,\Phi)$ of an odometer  $({\mathfrak{X}},\Gamma,\Phi)$ is well-defined, and that its type $\tau[\xi({\mathfrak{X}},\Gamma,\Phi)]$ is an invariant under return equivalence. The first claim is a special case of the proof of invariance under return equivalence, as explained later.

Suppose that the odometers $({\mathfrak{X}}, \Gamma, \Phi)$ and $({\mathfrak{X}'}, \Gamma', \Phi')$ are return equivalent. 
Then by Theorem~\ref{thm-algmodel}, the odometers $(X_{\infty}, \Gamma, \Phi_{\infty})$ and $(X_{\infty}', \Gamma', \Phi_{\infty}')$ are return equivalent as in Definition~\ref{def-return}.

Choose adapted sets $U \subset X_{\infty}$ and $U' \subset X'_{\infty}$ and a homeomorphism $h \colon U \to U'$ so that the induced homomorphism  
 $h_* \colon {\rm Homeo}(U) \to {\rm Homeo}(U')$ restricts to an isomorphism between the image   ${\mathcal H}_U = \Phi(\Gamma_U)$ with the image  ${\mathcal H}'_{U'} = \Phi'(\Gamma'_{U'})$.

Recall a standard construction of a  diagram of maps between adapted sets     \cite{FO2002},  \cite[Theorem~3.3]{DHL2016a}:
 \begin{align} \label{diag-lattice}
 \xymatrixcolsep{0pc}
\xymatrix{ 
X_{\infty} & \supset &    U\ar[d]^h  & \supset & \gamma_1 \cdot U_{\ell_1}\ar[d]^h & \supset & h^{-1}(\gamma'_1 \cdot U'_{\ell'_1})\ar[d]^h  & \supset & \gamma_{2} \cdot U_{\ell_2}\ar[d]^h  &  \supset &   h^{-1}(\gamma'_2 \cdot U'_{\ell'_2})\ar[d]^h  &  \supset &\cdots   \\
X'_{\infty} & \supset  & U'  & \supset  & h(\gamma_1 \cdot U_{\ell_1}) & \supset & \gamma'_1 \cdot U'_{\ell'_1}  &  \supset & h(\gamma_{2} \cdot U_{\ell_2})  &  \supset &   \gamma'_2 \cdot U'_{\ell'_2}  &   \supset &\cdots   } 
\end{align}
The sets $U_{\ell}$ are defined in \eqref{eq-openbasis}, and adapted to the action of $\Gamma$, and similarly for the sets $U'_{\ell'}$ which are  adapted to the action of $\Gamma'$. The subscripts and intermediate adapted sets are defined iteratively in the following.
We adopt the $``\cdot$'' notation for the actions, as it is clear from the context which action is being applied.
Recall that the action of $\Gamma_U$ on $U$ is minimal, as is the action of $\Gamma'_{U'}$ on $U'$. Denote $e_\infty = (e\Gamma_\ell) \in X_\infty$, and $e_{\infty}' = (e \Gamma_{\ell'}') \in X_\infty'$, where $\Gamma_\ell$ is the isotropy subgroup of $U_\ell$ and $\Gamma_\ell'$ is the isotropy subgroup of $U_{\ell'}'$, for each $\ell \geq 0$ and each $\ell' \geq 0$.

    We define the maps and indices in  Diagram~\eqref{diag-lattice} recursively. 
    
 The set $U$ is clopen, so there exists  $\gamma_1  \in \Gamma$ such that $\gamma_1 \cdot e_{\infty} \in U$.\\  Choose $\ell_1 > 0$ such that $\gamma_1 \cdot U_{\ell_1} \subset U$.
  
  The image $h(\gamma_1 \cdot U_{\ell_1}) \subset U'$ is clopen, so choose  $\gamma'_1 \in \Gamma'$ such that $\gamma'_1 \cdot e'_{\infty}  \in h(\gamma_1 \cdot U_{\ell_1})$. \\Choose $\ell'_1 > 0$ such that $\gamma'_1 \cdot U'_{\ell'_1} \subset  h(\gamma_1 \cdot U_{\ell_1})$.
  
  The image $h^{-1}(\gamma'_1 \cdot U'_{\ell'_1}) \subset \gamma_1 \cdot U_{\ell_1}$ is clopen, so choose $\gamma_2 \in \Gamma$ such that $\gamma_2 \cdot e_{\infty} \in   h^{-1}(\gamma'_1 \cdot U'_{\ell'_1}) \subset \gamma_1 \cdot U_{\ell_1}$. \\Choose $\ell_2 > \ell_1$ such that $\gamma_2 \cdot U_{\ell_2} \subset h^{-1}(\gamma'_1 \cdot U'_{\ell'_1})$.

  The image $h(\gamma_2 \cdot U_{\ell_2}) \subset \gamma'_1 \cdot U'_{\ell'_1}$ is clopen, , so choose $\gamma'_2 \in \Gamma'$ such that 
  $\gamma'_2 \cdot e'_{\infty} \in h(\gamma_2 \cdot U_{\ell_2})$. \\Choose $\ell'_2 > \ell'_1$ such that $\gamma'_2 \cdot U'_{\ell'_2} \subset h(\gamma_2 \cdot U_{\ell_2})$.

Continue this procedure recursively to obtain   Diagram~\eqref{diag-lattice} where  we have:

\begin{itemize}
\item increasing sequences $0< \ell_1 < \ell_2 < \ell_3 < \cdots$ and $0 < \ell'_1 < \ell'_2 < \ell'_3 < \cdots$,
\item a sequence $\{\gamma_{1},  \gamma_{2}, \gamma_{3}, \ldots\} \subset \Gamma$,
\item a sequence $\{\gamma'_{1},  \gamma'_{2}, \gamma'_{3}, \ldots\} \subset \Gamma'$ \ .
\end{itemize}
 
Observe that by choice we have $\gamma_{i+1} \cdot U_{\ell_{i+1}} \subset \gamma_{i} \cdot U_{\ell_i}$ and thus 
$\gamma_{\ell_i}^{-1} \gamma_{\ell_{i+1}} \in \Gamma_{\ell_i}$ for all $i \geq 1$. This recursion relation implies that the sequence $\{\gamma_i \mid i \geq 0\}$ converges in the profinite topology on $\Gamma$ induced by the subgroup chain $\{\Gamma_{\ell} \mid \ell \geq 0\}$, and similarly for $\{\gamma_i' \mid i \geq 0\}$ in the profinite topology on $\Gamma'$. If we denote the respective limits by $\overline{\{\gamma_i\}} \in \overline{\Gamma}$ and $\overline{\{\gamma_i'\}} \in \overline{\Gamma}'$, then we have
$$\overline{\{\gamma_i\}} \cdot e_{\infty} = h^{-1}(\overline{\{\gamma_i'\}} \cdot e'_{\infty}), \textrm{ and }\overline{\{\gamma'_i\}} \cdot e'_{\infty} = h(\overline{\{\gamma_i\}} \cdot e_{\infty}).$$

All of the sets appearing in Diagram~\eqref{diag-lattice} are adapted for their respective actions, so the set inclusions induce  corresponding subgroup chains of their isotropy groups, and these chains are interlaced:
\begin{align} \label{diag-latticegrps}
 \xymatrixcolsep{0pc}
\xymatrix{ 
\Gamma & \supset &  H=  \Gamma_U\ar[d]^{h_*}  & \supset & \gamma_1 \Gamma_{\ell_1} \gamma_1^{-1}\ar[d]^{h_*} & \supset & \Gamma_{h^{-1}(\gamma'_1 \cdot U'_{\ell'_1})}\ar[d]^{h_*}  & \supset & \gamma_2 \Gamma_{\ell_2} \gamma_2^{-1}\ar[d]^{h_*}  &  \supset &   \Gamma_{h^{-1}(\gamma'_2 \cdot U'_{\ell'_2})}\ar[d]^{h_*}  &  \supset &\cdots   \\
\Gamma' & \supset  & H' = \Gamma'_{U'}  & \supset  & \Gamma'_{h(\gamma_1 \cdot U_{\ell_1})} & \supset & \gamma'_1 \Gamma'_{\ell'_1} (\gamma'_1)^{-1}  &  \supset & \Gamma'_{h(\gamma_2 \cdot U_{\ell_2})}  &  \supset &   \gamma'_2 \Gamma'_{\ell'_2}(\gamma'_2)^{-1}  &   \supset &\cdots   } 
\end{align}

Conjugation does not change the index of a subgroup, so $[\Gamma : \Gamma_{\ell_j}] = [\Gamma : \gamma_j \Gamma_{\ell_j} \gamma_j^{-1}]$ for all $j \geq 1$, and likewise we have
$[\Gamma' : \Gamma'_{\ell'_j}] = [\Gamma' : \gamma'_j \Gamma'_{\ell'_j} (\gamma'_j)^{-1}]$.
 It then follows  that 
 \begin{equation}    \label{eq-factor1} 
 [\Gamma :  \Gamma_{\ell_j}]     =    [\Gamma :   \Gamma_{\ell_1}  ] [\Gamma_{\ell_1}: \Gamma_{\ell_j}]   =   [\Gamma :   \Gamma_{\ell_1} ][\gamma_1 \Gamma_{\ell_1} \gamma_1^{-1} : \gamma_j \Gamma_{\ell_j} \gamma_j^{-1}]     \ ,    
  \end{equation}
  \begin{equation}       \label{eq-factor2}
  [\Gamma' :  \Gamma'_{\ell'_j}]   =      [ \Gamma' :  \Gamma'_{\ell'_1}  ]   [\Gamma'_{\ell'_1} :  \Gamma'_{\ell'_j}  ]  
  =    [ \Gamma'  :   \Gamma'_{\ell_1}][\gamma'_1 \Gamma'_{\ell'_1} (\gamma'_1)^{-1} : \gamma'_j \Gamma'_{\ell'_j} (\gamma'_j)^{-1}]    \ .  
 \end{equation}
   We now show the key fact for the proof of Theorem~\ref{thm-main2}.
   \begin{lemma}\label{lem-key}
   For all $j > i > 0$, 
   \begin{equation}\label{eq-indices}
  [\Gamma_{\ell_i}: \Gamma_{\ell_j}]   =   [\Gamma'_{h(\gamma_i \cdot U_{\ell_i})} : \Gamma'_{h(\gamma_j \cdot U_{\ell_j})}].
\end{equation}
     \end{lemma}
   \proof
 The isotropy group $\Gamma_{\ell_i}$ acts minimally on the clopen set $U_{\ell_i}$, and its action translates the clopen subset 
 $U_{\ell_j} \subset U_{\ell_i}$ to give a partition of this set.
   The index $[\Gamma_{\ell_i}: \Gamma_{\ell_j}]$ equals the number of clopen subsets in this partition.
   It follows that the action of $\gamma_i \Gamma_{\ell_i} \gamma_i^{-1}$ on $\gamma_i \cdot U_{\ell_i}$ partitions this set into 
   $[\Gamma_{\ell_i}: \Gamma_{\ell_j}]$ translates of the clopen subset $\gamma_j \cdot U_{\ell_j}$.

   The homeomorphism $h \colon U \to U'$ restricts to a homeomorphism $h \colon \gamma_i \cdot U_{\ell_i} \to h(\gamma_i \cdot U_{\ell_i})$ which induces a conjugacy of the action of $\Gamma_{\gamma_i \cdot U_{\ell_i}}$ with the action of 
   $\Gamma'_{h(\gamma_i \cdot U_{\ell_i})}$ on $ h(\gamma_i \cdot U_{\ell_i}) \subset U'$. Thus, the translates of the clopen subset 
 $h(\gamma_j \cdot U_{\ell_j})$   partition $h(\gamma_i \cdot U_{\ell_i})$ into $[\Gamma_{\ell_i}: \Gamma_{\ell_j}]$ clopen subsets, and so the identity \eqref{eq-indices} follows.
   \endproof
   
 We now prove that  $\tau[\xi(X_{\infty}, \Gamma, \Phi_{\infty})] = \tau[\xi(X_{\infty}', \Gamma', \Phi_{\infty}')]$. Recall that
 \begin{equation}\label{eq-LCM1}
\xi({\mathcal G}_{\mathcal U}) = {\rm lcm} \{ [\Gamma : \Gamma_{\ell}] \mid  \ell > 0\}   = {\rm lcm} \{ [\Gamma : \Gamma_{\ell_j}] \mid j > 0 \} \ ,
\end{equation}
 \begin{equation}\label{eq-LCM2}
\xi({\mathcal G}'_{{\mathcal U}'}) = {\rm lcm} \{ [\Gamma' : \Gamma'_{\ell}] \mid  \ell > 0\}  = {\rm lcm} \{ [\Gamma' : \Gamma'_{\ell'_j}] \mid j > 0 \} \ .
\end{equation}

 Then by     \eqref{eq-indices},  for $j > 1$ we have
 \begin{equation}\label{eq-changeofform}
[ \Gamma : \Gamma_{\ell_j} ]  = [ \Gamma : \Gamma_{\ell_1}] \ [\Gamma_{\ell_1} : \Gamma_{\ell_j}]  = [ \Gamma : \Gamma_{\ell_1}] \  [\Gamma'_{h(\gamma_1 \cdot U_{\ell_1})} : \Gamma'_{h(\gamma_j \cdot U_{\ell_j})}]  \ .
\end{equation}
 Then calculate,  recalling the inclusions $U' \supset h(\gamma_1 \cdot U_{\ell_1}) \supset \gamma_1' \cdot U_{\ell_1}' \supset h(\gamma_j \cdot U_{\ell_j}) \supset \gamma_j' \cdot U_{\ell_j'}'$ from \eqref{diag-lattice}:
 \begin{eqnarray}\label{eq-bigmess}
\lefteqn{ [\Gamma' : \Gamma'_{\ell'_j}]   =    [\Gamma' : \gamma_j'\Gamma'_{\ell'_j} (\gamma_j')^{-1}] } \\
& = &  [\Gamma' :  \Gamma'_{h(\gamma_1 \cdot U_{\ell_1})} ]  \ 
 [\Gamma'_{h(\gamma_1 \cdot U_{\ell_1})} :   \Gamma'_{\gamma'_1 \cdot U'_{\ell'_1}} ] \ [\Gamma'_{\gamma'_1 \cdot U'_{\ell'_1}} : \Gamma'_{h(\gamma_{j} \cdot U_{\ell_{j}})}]  \  [ \Gamma'_{h(\gamma_{j} \cdot U_{\ell_{j}})}  : \Gamma'_{\gamma'_j \cdot U'_{\ell'_j}}] \nonumber \\
& = &  [\Gamma' :  \Gamma'_{h(\gamma_1 \cdot U_{\ell_1})} ]  \ 
 [\Gamma'_{h(\gamma_1 \cdot U_{\ell_1})}     : \Gamma'_{h(\gamma_{j} \cdot U_{\ell_{j}})}] \ [\Gamma'_{h(\gamma_{j} \cdot U_{\ell_{j}})} : \Gamma'_{\gamma'_j \cdot U'_{\ell'_j}}] .\nonumber
\end{eqnarray}
Then by \eqref{eq-changeofform},     
the last line of \eqref{eq-bigmess} yields
 \begin{equation} \label{eq-balance}  
 [\Gamma : \Gamma_{\ell_1}] \  [\Gamma' : \Gamma'_{\ell'_{j}}]  =  [\Gamma' :  \Gamma'_{h(\gamma_1 \cdot U_{\ell_1})} ]  \ [ \Gamma : \Gamma_{\ell_j} ] \ 
 [\Gamma'_{h(\gamma_{j} \cdot U_{\ell_{j}})} : \Gamma'_{\gamma'_{j} \cdot U'_{\ell'_{j}}}].
     \end{equation}
  
Thus, the index   $[ \Gamma : \Gamma_{\ell_j}]$ divides  $[\Gamma : \Gamma_{\ell_1}] \  [\Gamma' : \Gamma'_{\ell'_{j}}]$ for all $j > 1$. In particular, 
$[ \Gamma : \Gamma_{\ell_j}]$ divides $[\Gamma : \Gamma_{\ell_1}] \ \xi({\mathcal P}')$ for all $j > 1$, and so 
  $\xi({\mathcal P})$ divides $[\Gamma : \Gamma_{\ell_1}] \ \xi({\mathcal P}')$, hence $\tau[{\mathcal P}] \leq \tau[{\mathcal P}']$. 

Next, repeat these calculations for $j > 1$, starting with 
 \begin{equation}\label{eq-changeofform2}
[ \Gamma : \Gamma_{\ell_{j+1}} ]  =   [ \Gamma : \Gamma_{h^{-1}(\gamma'_1 \cdot U_{\ell'_1})}] \  [\Gamma_{h^{-1}(\gamma'_1 \cdot U_{\ell'_1})} : \Gamma_{h^{-1}(\gamma'_j \cdot U'_{\ell'_j})}] \ [\Gamma_{h^{-1}(\gamma'_j \cdot U'_{\ell'_j})} :  \Gamma_{\gamma_{j+1} \cdot U_{\ell_{j+1}}}]   \ .
\end{equation}
Lemma~\ref{lem-key} can be applied to the inverse map $h^{-1} \colon U' \to U$ as well, to obtain that
for all $j > i > 0$, 
   \begin{equation}\label{eq-indicesinv2}
  [\Gamma'_{\ell'_i}: \Gamma'_{\ell'_j}]   =   [\Gamma_{h^{-1}(\gamma'_i \cdot U_{\ell'_i})} : \Gamma_{h^{-1}(\gamma'_j \cdot U'_{\ell'_j})}] 
\end{equation}
and so
\begin{equation}\label{eq-changeofforminv2}
[ \Gamma : \Gamma_{\ell_{j+1}} ]  =   [ \Gamma : \Gamma_{h^{-1}(\gamma'_1 \cdot U_{\ell'_1})}] \   [\Gamma'_{\ell'_1}: \Gamma'_{\ell'_j}]   \ [\Gamma_{h^{-1}(\gamma'_j \cdot U'_{\ell'_j})} :  \Gamma_{\gamma_{j+1} \cdot U_{\ell_{j+1}}}]   \ .
\end{equation}

Thus, the index    $[\Gamma'_{\ell'_1}: \Gamma'_{\ell'_j}]$  divides  $[\Gamma : \Gamma_{\ell_{j+1}}]$  for all $j > 1$ and so $[\Gamma'_{\ell'_1}: \Gamma'_{\ell'_j}]$ divides $\xi({\mathcal G}_{\mathcal U})$. It follows that 
  $\xi({\mathcal G}'_{{\mathcal U}'})$ divides $\xi({\mathcal G}_{\mathcal U})$, hence $\tau[{\mathcal G}'_{{\mathcal U}'}] \leq \tau[{\mathcal G}_{\mathcal U}]$. 
Reversing the roles of the group chains ${\mathcal G}_{\mathcal U}$ and ${\mathcal G}'_{{\mathcal U}'}$, we obtain that $\tau[{\mathcal G}_{\mathcal U}] =  \tau[{\mathcal G}'_{{\mathcal U}'}]$, and thus $\tau[{\mathfrak{X}}, \Gamma, \Phi)] = \tau[{\mathfrak{X}'}, \Gamma', \Phi']$.  This completes the proof that the types of   return equivalent actions are equal.

Finally, to complete the proof of Theorem~\ref{thm-main2},  we show that $\xi({\mathcal G}_{\mathcal U})=  \xi({\mathcal G}'_{{\mathcal U}'})$ for     group chains ${\mathcal G}_{\mathcal U}$ and ${\mathcal G}'_{{\mathcal U}'}$ associated with  $({\mathfrak{X}},\Gamma,\Phi)$.  We proceed as in the above proof, and note that  while $\Gamma = \Gamma'$,   the group chains define distinct  inverse limit spaces $X_{\infty}$ and $X'_{\infty}$. In the above calculations, take $U= X_{\infty}$ and $U'=X'_{\infty}$, $\gamma_1 = e \in \Gamma$ with $\ell_1 = 0$, and $\gamma'_1 = e' \in \Gamma'$ with $\ell'_1 = 0$. Then the terms $ [\Gamma : \Gamma_{\ell_1}] =1$ and  $[\Gamma' : \Gamma'_{\ell'_1}] = 1$. 
Then by \eqref{eq-balance} we have $[ \Gamma : \Gamma_{\ell_j} ]$ divides $[\Gamma' : \Gamma'_{\ell'_{j}}]$, and by \eqref{eq-changeofforminv2} we have  $[\Gamma' : \Gamma'_{\ell'_j}]$ divides $[ \Gamma : \Gamma_{\ell_{j+1}} ]$. It follows that 
$\xi({\mathcal G}_{\mathcal U})=  \xi({\mathcal G}'_{{\mathcal U}'})$. That is,  the Steinitz  order is independent of the choice of a group associated to the action.

  \subsection{Proof of Theorem~\ref{thm-main2b}}\label{subsec-main2b}
  
 We must show that the typesets of odometers are well-defined, and invariant under return equivalence modulo the equivalence relation on typesets in Definition~\ref{def-commensurable}. As above, we show the result for return equivalence first, then deduce the result for conjugation from the proof of this. The proof begins   as in Section~\ref{subsec-proof-main2}.

Suppose that the odometers $({\mathfrak{X}}, \Gamma, \Phi)$ and $({\mathfrak{X}'}, \Gamma', \Phi')$ are return equivalent, 
then the odometers $(X_{\infty}, \Gamma, \Phi_{\infty})$ and $(X_{\infty}', \Gamma', \Phi_{\infty}')$ are return equivalent.
Choose adapted sets $U \subset X_{\infty}$ and $U' \subset X'_{\infty}$ and a homeomorphism $h \colon U \to U'$ so that the induced homomorphism  
 $h_* \colon {\rm Homeo}(U) \to {\rm Homeo}(U')$ restricts to an isomorphism between the image   ${\mathcal H}_U = \Phi(\Gamma_U)$ with the image  ${\mathcal H}'_{U'} = \Phi'(\Gamma'_{U'})$
that induces an isomorphism between the actions $(U, {\mathcal H}_{U}, \Phi_{U})$ and $(U', {\mathcal H}'_{U'}, \Phi'_{U'})$.

Next, choose basepoints $x \in U$ and $x' \in U'$. Choose an   adapted neighborhood basis ${\mathcal U} = \{U_{\ell} \subset {\mathfrak{X}}  \mid \ell > 0\}$ at $x$ for the action $\Phi$ as in Definition~\ref{def-adaptednbhds}, and an     adapted neighborhood basis ${\mathcal U}' = \{U'_{\ell} \subset {\mathfrak{X}}'  \mid \ell > 0\}$ at $x'$ for the action $\Phi'$. We can assume that $U_1 \subset U$ and $U'_1 \subset U'$.

 Form the  group chain $\displaystyle  {\mathcal G} = \{\Gamma_0 \supset \Gamma_1 \supset \cdots \}$, where 
$\Gamma_{\ell} = \Gamma_{U_{\ell}}$   with $\Gamma_0 = \Gamma$, and   the group chain  $\displaystyle  {\mathcal G}' = \{\Gamma'_0 \supset \Gamma'_1 \supset \cdots \}$, where 
$\Gamma'_{\ell} = \Gamma'_{U'_{\ell}}$   with $\Gamma'_0 = \Gamma'$.

Set $H = \Gamma_U$ and $H' = \Gamma'_{U'}$.   We   will show that the $H$-restricted typeset     
 $\Xi_H[X_{\infty}, \Gamma, \Phi]$ and the $H'$-restricted typeset   $\Xi_{H'}[X'_{\infty} , \Gamma' , \Phi']$ are equal.

Given $\gamma \in \Gamma$, for  any positive integer $m > 0$ the group $\langle \gamma^m \rangle$ is a subgroup of finite index in $\langle \gamma \rangle$.  Thus  $\xi^H(\gamma) = m \xi^H(\gamma^m)$, and $\tau[\xi^H(\gamma)] =   \tau[\xi^H(\gamma^m)]$.
 Since $H$ has finite index in $\Gamma$, for any $\gamma \in \Gamma$ there exists $m > 0$ such that $\gamma^m \in H$. Thus we have
$\Xi_H[{\mathcal G}] = \Xi_H[H \cap {\mathcal G}]$, so it suffices to consider $\gamma \in H$, and likewise for $\gamma' \in H'$.

Let   sequences $\{\gamma_i \mid i \geq 1\}$, $\{\gamma'_i \mid i \geq 1\}$, $\{\ell_i \mid i \geq 1\}$ and $\{\ell'_i \mid i \geq 1\}$ be chosen as in Section~\ref{subsec-proof-main2}, 
resulting in the diagram \eqref{diag-lattice} of  adapted sets, and   the diagram \eqref{diag-latticegrps} of group inclusions.  Note that as we assume $e_{\infty} \in U$ we can choose $\gamma_1$ to be the identity, and likewise as $e'_{\infty} \in U'$ we choose $\gamma'_1$ to be the identity. By choice we have   $\gamma_{i+1} \cdot U_{\ell_{i+1}} \subset \gamma_{i} \cdot U_{\ell_i}$ and thus 
$\gamma_{\ell_i}^{-1} \gamma_{\ell_{i+1}} \in \Gamma_{\ell_i}$ for all $i \geq 1$. In particular, this implies that $\gamma_i \in H$ for all $i \geq 1$. The analogous conclusion holds, that $\gamma'_i \in H'$ for all $i \geq 1$.

For notational convenience, set $H_i = \Gamma_{\gamma_i \cdot U_{\ell_i}} = \gamma_i \Gamma_{\ell_i} \gamma_i^{-1}$ and $H'_{i} = \Gamma'_{h(\gamma_i \cdot U_{\ell_i})}$ for $i > 0$. Then $H_{j} \subset H_i \subset H = \Gamma_U$ and 
$H'_{j} \subset H'_i \subset H' = \Gamma'_{U'}$ for $j > i > 0$.

  Let $\{X_{\infty} : U\}$ denote the set of translates of  the adapted set $U$. Then for $H = \Gamma_U$ we have the set equality 
  $\{X_{\infty} : U\} = \Gamma/H$. 
  Set  $m = \#(\Gamma/H)!$ which       is the order of the group of permutations on the set $\{X_{\infty} : U\}$. 
 Then for  $\gamma \in \Gamma$  the  action of $\gamma^m$ on the set $\{X_{\infty} : U\}$ is the identity.  This implies that $\gamma^m \in C^{\Gamma}_H$ the core of $H$.

Given a group chain ${\mathcal G}$ and   $H \subset \Gamma$ a subgroup with $\Gamma_{\ell} \subset H$ for some $\ell \geq 1$, the \emph{$H$-restricted order} of $\gamma \in \Gamma$ is the  Steinitz order with respect to the chain ${\mathcal C}_{\gamma}^H$, 
\begin{equation}\label{eq-Htypegamma2}
\xi^H(\gamma) = {\rm lcm} \{   \#( \langle \gamma \rangle /C_{\gamma, \ell}^H  ) \mid \ell > 0 \}  \ .
\end{equation}
The $H-$\emph{restricted typeset} for the chain ${\mathcal G}$  is the collection
\begin{equation}\label{eq-Htypeaction2}
\Xi_H[{\mathcal G}] = \{ \tau[\xi^H(\gamma)] \mid \gamma \in \Gamma\} \ .
\end{equation}

 \begin{lemma}
 $\tau[\xi^H(\gamma^m)]=     \tau[\xi^H(\gamma)]$
 \end{lemma}
 \proof
Let the cyclic group $\langle \gamma \rangle /C_{\gamma, \ell}^H$ have order $m_{\ell}$. Then $\gamma^m$ generates a subgroup whose index is the greatest common divisor $n_{\ell} = {\rm gcd}(m,m_{\ell})$. We have $m_{\ell} | m_{\ell +1}$ hence $n_{\ell} \leq n_{\ell + 1}$, and also that   $n_{\ell} \leq m$. Thus,  the sequence $\{n_{\ell} \mid \ell \geq 1\}$ has an upper bound $n_{\infty} \leq m$ which is realized for $\ell$ sufficiently large. It follows that 
 $n_{\infty} \cdot \xi^H(\gamma^m)=     \xi^H(\gamma)$,  hence  $\tau[\xi^H(\gamma^m)]=     \tau[\xi^H(\gamma)]$.
\endproof

So without loss of generality it suffices to consider $\gamma \in C^{\Gamma}_H \subset H$.

The homeomorphism $h \colon U \to U'$ induces an isomorphism $h_* \colon {\mathcal H}_U \to {\mathcal H}'_{U'}$.
Thus there exists $\gamma' \in H'$ whose action $\Phi'_{U'}(\gamma')$ on $U'$ equals the image $\phi'_{\gamma} = h_*(\Phi_U(\gamma))$. 
 We show that the $H'$-restricted type of $\gamma'$ equals the $H$-restricted type of $\gamma$.

Let  $V \subset W \subset X_{\infty}$ be adapted subsets. The restricted action of $\Gamma_W$ on $W$  is minimal, hence the translates of $V$ define a clopen partition of $W$. Let $\{W:V\}$ denote the set of elements in this partition, and let $|W:V| = \#\{W:V\}$ denote the cardinality of the set of translates.

 The action $\Phi$ induces a map  $\Phi_{V}^W \colon \Gamma_W \to {\rm Aut}(\{W:V\})$ into the permutations of the set $\{W:V\}$ of translates of $V$ in $W$. The kernel of the map $\Phi_{V}^W$ is   denoted by $C_V^W \subset \Gamma_V$ and equals the  normal core of $\Gamma_V$ as a subgroup of  $\Gamma_W$. Thus the index $[\Gamma_W: \Gamma_V] = |W:V|$, and  $[\Gamma_W: C^W_V] = \# {\rm Image}(\Phi_{V}^W(\Gamma_W))$.

 Now apply this observation to the action of $H = \Gamma_U$ on $U$. For  each $j \geq 1$,  the action $\Phi$  induces a map $\widehat{\Phi}_{U_{\ell_j}}^U \colon H \to {\rm Aut}(\{U : \gamma_j \cdot U_{\ell_j}\})$ which permutes the elements of this partition. The kernel of  the action map $\widehat{\Phi}_{U_{\ell_j}}^U$ is   the  subgroup 
 \begin{equation}
C^H_j =   \bigcap_{\delta \in H} ~  \delta  \Gamma_{\gamma_j \cdot U_{\ell_j}}  \delta^{-1} =   \bigcap_{\delta \in H} ~( \delta \gamma_j) \Gamma_{\ell_j} (\delta \gamma_j )^{-1}    =   \bigcap_{\delta \in H} ~ \delta  \Gamma_{\ell_j} \delta^{-1}   \ .
\end{equation}

 Set $C_{\gamma, j}^H = \langle \gamma \rangle \cap C_{j}^H$, then the subgroup 
  $\langle \gamma \rangle /C_{\gamma, j}^H  \subset H/C_{j}^H$ is mapped injectively into ${\rm Aut}(\{U : U_{\ell_j}\})$.
Thus, $\#( \langle \gamma \rangle /C_{\gamma, j}^H  ) = \# {\rm Image}(\Phi_{U_{\ell_j}}^U(\langle \gamma \rangle))$,  and so we have
  \begin{equation}\label{eq-counting}
\xi^H(\gamma) = {\rm lcm} \{   \#( \langle \gamma \rangle /C_{\gamma, j}^H  ) \mid j > 0 \}  = 
{\rm lcm} \{  \# {\rm Image}(\widehat{\Phi}_{U_{\ell_j}}^U(\langle \gamma \rangle)) \mid j > 0 \} \ .
\end{equation}

We next use the conjugation by $h$ to relate the order of a subgroup of ${\rm Aut}(\{U : V\})$ with the order of a subgroup of ${\rm Aut}(\{U' : V'\})$ for appropriate choices of adapted subsets $V$ and $V'$.     This is analogous to the idea behind the proof of Lemma~\ref{lem-key}.

From  \eqref{diag-lattice}, for $j \geq 1$ we have
\begin{equation}\label{eq-containment}
 U'    \supset    h(\gamma_1 \cdot U_{\ell_1})   \supset   \gamma'_1 \cdot U'_{\ell'_1}        \supset   \cdots 
    \supset      h(\gamma_{j} \cdot U_{\ell_j})  \supset    \gamma'_j \cdot U'_{\ell'_j}  \supset   \cdots   \ .
\end{equation}

Using that  $\gamma'_j \cdot U'_{\ell'_j} \subset h(\gamma_{j} \cdot U_{\ell_j})$  we have 
\begin{equation}
\gamma'_{j} \Gamma'_{\ell_{j}}(\gamma'_{j})^{-1}  = \Gamma'_{\gamma'_j \cdot U'_{\ell'_j}} \subset \Gamma'_{h(U_{\ell_j})} = H'_j \subset H' \ , 
\end{equation}
  and so $C^{H'}_j = C^{H'}_{\gamma'_{j} \Gamma'_{\ell_{j}}(\gamma'_{j})^{-1}} \subset C^{H'}_{H'_j}$. Then
  \begin{equation}\label{eq-chainintersections}
   \langle \gamma' \rangle  \supset \Big( \langle \gamma' \rangle \cap C^{H'}_{H'_j} \Big) \supset \Big( \langle \gamma' \rangle \cap C^{H'}_{j} \Big) = C^{H'}_{\gamma',j}.
  \end{equation}
 
The action of $\Phi'_{U'}(\gamma')$ on $U'$   translates  the clopen set $h(\gamma_j \cdot U_{\ell_j})$ within the clopen set $U'$, thus
\begin{equation}\label{eq-counting2ell}
\#( \langle \gamma \rangle /C_{\gamma, j}^H  )   = 
 \# {\rm Image}(\widehat{\Phi}_{U_{\ell_j}}^U(\langle \gamma \rangle)) = 
 \# {\rm Image}(\widehat{\Phi}_{h(U_{\ell_j})}^{U'}(\langle \gamma' \rangle)) = 
\#( \langle \gamma' \rangle /(\langle \gamma' \rangle \cap C^{H'}_{H'_j}  )  )  \ . 
\end{equation}

Using that  $\gamma'_j \cdot U'_{\ell'_j} \subset h(\gamma_{j} \cdot U_{\ell_j})$  we have 
\begin{equation}
\gamma'_{j} \Gamma'_{\ell_{j}}(\gamma'_{j})^{-1}  = \Gamma'_{\gamma'_j \cdot U'_{\ell'_j}} \subset \Gamma'_{h(U_{\ell_j})} = H'_j \subset H' \ , 
\end{equation}
  and so $C^{H'}_j = C^{H'}_{\gamma'_{j} \Gamma'_{\ell_{j}}(\gamma'_{j})^{-1}} \subset C^{H'}_{H'_j}$.
  We then have, using \eqref{eq-chainintersections},
  \begin{eqnarray}
 \#( \langle \gamma' \rangle /C_{\gamma', j}^{H'}  ) & = &  \#( \langle \gamma' \rangle /(\langle \gamma' \rangle \cap C^{H'}_{H'_j}  )  ) \cdot \#((\langle \gamma' \rangle \cap C^{H'}_{H'_j}  ) / C_{\gamma', j}^{H'}  ) \label{eq-divided} \\
 & = &  \#( \langle \gamma \rangle /C_{\gamma, j}^H  ) \cdot  \#((\langle \gamma' \rangle \cap C^{H'}_{H'_j}  ) / (\langle \gamma' \rangle \cap C^{H'}_{j}  ) ) \nonumber \ .
\end{eqnarray}

Thus, we have that $\#( \langle \gamma \rangle /C_{\gamma, j}^H  )$ divides $\#( \langle \gamma' \rangle /C_{\gamma', j}^{H'}  ) $ for all $j \geq 1$. It follows that $\tau_H[\gamma] = \tau[\xi^{H}(\gamma)] \leq \tau[\xi^{H'}(\gamma')]  = \tau_{H'}[\gamma']$. Then by reversing these calculations, starting with $\gamma'$ chosen as above, we obtain $\tau_{H'}[\gamma'] \leq \tau_H[\gamma]$ and hence $\tau_H[\gamma] = \tau_{H'}[\gamma']$.
 
 We have thus shown that for each $\gamma \in \Gamma$ there exists $\gamma' \in \Gamma'$ with the same restricted type. Reversing this process, we deduce that $\Xi_H[X_{\infty},\Gamma,\Phi] = \Xi_{H'}[X'_{\infty},\Gamma',\Phi']$. 
 It follows that $\Xi[{\mathcal G}]$ and $\Xi[{\mathcal G}']$ are commensurable as in
 Definition~\ref{def-commensurable}, which proves Theorem~\ref{thm-main2b}.

    \subsection{Proof of Theorem~\ref{thm-main2a}}

  The claim of Theorem~\ref{thm-main2a} is that for each   $\gamma \in \Gamma$, there is a well-defined type $\tau[\gamma]$. In fact, we show the stronger result, that the Steinitz order $\xi(\gamma)$    is independent of the choice of an adapted basis ${\mathcal U}$ for the odometer $({\mathfrak{X}}, \Gamma, \Phi)$. Given two adapted bases ${\mathcal U}$ and ${\mathcal U}'$,  the algebraic models they define, $(X_{\infty}, \Gamma, \Phi_{\infty})$ and $(X'_{\infty}, \Gamma, \Phi_{\infty}')$, are isomorphic. 

       In the above proof of Theorem~\ref{thm-main2b}, we let $U = X_{\infty}$ and $U' = X'_{\infty}$ and so $H = \Gamma_U = \Gamma$ and $H' = \Gamma'_{U'} = \Gamma'$. Then observe that for $\ell > 0$ the normal subgroups $C^H_{\ell} = C_{\ell}$ and $C^{H'}_{\ell} = C'_{\ell}$. It follows that for $\gamma \in \Gamma$ and $\gamma' \in \Gamma'$ chosen as in the proof of part (3b) we have $\xi^H(\gamma) = \xi(\gamma)$ and $\xi^{H'}(\gamma') = \xi(\gamma')$, and the proof shows that $\xi(\gamma) = \xi(\gamma')$. It follows that the set of orders $\{ \xi(\gamma) \mid \gamma \in \Gamma\}$ is independent of the choice of an algebraic model for the action, and so is an invariant of the isomorphism class of the action $({\mathfrak{X}}, \Gamma, \Phi)$, as was to be shown.

\subsection{Proof of Corollary~\ref{cor-abelian}}
Suppose that  both $\Gamma$ and $\Gamma'$ are abelian, and the actions are return equivalent. Then every group $\Gamma_{\ell}$ in a subgroup chain ${\mathcal G}_{{\mathcal U}}$ in $\Gamma$ is normal. Then for any subgroup $H \subset \Gamma$, the normal cores satisfy $C^H_{\ell} = C_{\ell} = \Gamma_{\ell}$. Thus, in \eqref{eq-counting}, the restricted type $\xi^H(\gamma) = \xi(\gamma)$. Similarly, we have $\xi^{H'}(\gamma') = \xi(\gamma')$ for $\gamma' \in \Gamma'$. Thus, the above proof shows in this case that the typesets $\Xi[{\mathfrak{X}},\Gamma,\Phi]$ and $\Xi[{\mathfrak{X}}',\Gamma',\Phi']$ are equal, which is the claim of  Corollary~\ref{cor-abelian}.

It is surprising perhaps, that the conclusion of Corollary~\ref{cor-abelian} need not hold if one of the groups is virtually abelian, but not abelian, as illustrated in Example~\ref{ex-twistedabelian}. 
The idea of Example~\ref{ex-twistedabelian} is that we add to an abelian group a single element that normalizes it, does not commute with it, and destroys the equality $C^H_{\ell} = C_{\ell}$. Then the types of elements are no longer equal to their  $H$-restricted types.

 \subsection{Proof of Theorem~\ref{thm-renorm}}
Let  $\varphi \colon \Gamma \to \Gamma$ be a self-embedding with finite index image.  Define a group chain ${\mathcal G}_{\varphi}$ by setting $\Gamma_0 = \Gamma$ and then inductively defining $\Gamma_{\ell +1} = \varphi(\Gamma_{\ell}) \subset \Gamma_{\ell}$ for $\ell \geq 0$. 

Then with   $H = \Gamma_{k}$ for   $k > 0$ and  $\ell > k$,
\begin{equation}
C^H_{\ell} = \bigcap_{\delta \in H} \ \delta^{-1} \Gamma_{\ell} \delta = \bigcap_{\delta \in \varphi^{k}(\Gamma)} \ \delta^{-1} \varphi^{\ell}(\Gamma) \delta = \varphi^k\left( \bigcap_{\delta \in \Gamma} \ \delta^{-1} \varphi^{\ell-k}(\Gamma)  \delta \right) = \varphi^k(C_{\ell-k}) \ ,  
\end{equation}

\bigskip
\bigskip

and so $\Gamma_{\ell}/C^H_{\ell} = \varphi^k(\Gamma_{\ell-k}/C_{\ell-k})$. Then for $\gamma \in H = \Gamma_{k}$ set $\overline{\gamma} = \varphi^{-k}(\gamma)$, then have
\begin{equation}
\xi^H(\gamma) = {\rm lcm} \{   \#( \langle \gamma \rangle /C_{\gamma, \ell}^H  ) \mid \ell > k \}  = 
{\rm lcm} \{   \#( \langle \overline{\gamma} \rangle /C_{\overline{\gamma}, \ell-k}  ) \mid \ell -k > 0 \} = \xi(\overline{\gamma}) \ .
\end{equation}
Then as $\Gamma = \varphi^{-k}(H)$, we have  $\Xi_H[X_{\varphi}, \Gamma, \Phi_{\varphi}] = \Xi[X_{\varphi}, \Gamma, \Phi_{\varphi}]$. 
  
Let $(X_{\varphi}, \Gamma, \Phi_{\varphi})$ and $(X_{\varphi'}, \Gamma', \Phi_{\varphi'})$ be odometers associated to renormalizations $\varphi \colon \Gamma \to \Gamma$ and $\varphi' \colon \Gamma' \to \Gamma'$. Assume the   actions  are return equivalent by a homeomorphism $h \colon U \to U'$. Then  choose $\delta \in \Gamma$ with $\delta \cdot e_{\infty} \in U$ and $\delta' \in \Gamma'$ with $\delta' \cdot e'_{\infty} \in U'$.
Then $\varphi^{\delta} = \Phi(\delta) \circ \varphi \circ \Phi(\delta^{-1})$ is a renormalization of $\Gamma$ with group chain  $\Gamma^{\delta}_{\ell} = \delta \Gamma_{\ell} \delta^{-1}$. The group $\Gamma^{\delta}_{\ell}$ stabilizes the translate $\delta \cdot U_{\ell}$. As $\delta \cdot e_{\infty} \in U$, there exists $k > 0$ such that $\delta \cdot U_k \subset U$. Repeat this argument for the renormalization $\varphi'$ to obtain a $\delta' \in \Gamma'$ and $k ' > 0$ such that $\delta' \cdot U'_{k'} \subset U'$.

Then proceed as in the proof in Section~\ref{subsec-main2b} with $H = \delta \Gamma_{k} \delta^{-1}$ and $H' = \delta' \Gamma'_{k'} (\delta')^{-1}$ to obtain that   
 $$\Xi[X_{\varphi}, \Gamma, \Phi_{\varphi}] = \Xi_H[X_{\varphi}, \Gamma, \Phi_{\varphi}] =  \Xi_{H'}[X'_{\varphi'}, \Gamma', \Phi_{\varphi'}] =  \Xi[X'_{\varphi'}, \Gamma', \Phi_{\varphi'}],$$ 
 as claimed in Theorem~\ref{thm-renorm}. In particular, this implies that the typeset $\Xi[X_{\varphi}, \Gamma, \Phi_{\varphi}]$ is an invariant of the isomorphism class of the renormalizable Cantor action $(X_{\varphi}, \Gamma, \Phi_{\varphi})$.

 \section{Basic Examples}\label{sec-abelian}

 \subsection{Virtually abelian actions} 
 A group $\Gamma$ is \emph{virtually abelian} if it admits a finite-index subgroup $A \subset \Gamma$  which is abelian. It is straightforward to construct ${\mathbb Z}^n$-odometer actions with prescribed spectra, and yet they illustrate several basic properties of the type and typeset invariants.  
 
 \begin{ex}\label{ex-toral}
 {\rm 
Consider the case $n=1$. Choose two disjoint sets of distinct primes,
$$\pi_f = \{q_1 , q_2, \ldots \} \quad , \quad \pi_{\infty} = \{p_1 , p_2, \ldots\}$$
where $\pi_f$ and $\pi_{\infty}$ can be chosen to be finite or infinite sets, and either $\pi_f$ is infinite, or $\pi_{\infty}$ is non-empty. Choose  multiplicities $n(q_i) \geq 1$ for the primes in $\pi_f$.  For each $\ell > 0$, define a subgroup of $\Gamma = {\mathbb Z}$ by 
\begin{equation}\label{eq-lotsaprimes}
\Gamma_{\ell} = \{q_1^{n(q_1)} q_2^{n(q_2)} \cdots q_{\ell}^{n(q_{\ell})} \cdot p_1^{\ell} p_2^{\ell} \cdots p_{\ell}^{\ell} \cdot n \mid n \in {\mathbb Z} \} \ .
\end{equation}
If $\pi_{\infty}$ is a finite set, then we use the convention that $p_{\ell} =1$ in \eqref{eq-lotsaprimes} when $p_{\ell}$ is not defined by the listing of $\pi_{\infty}$.
The completion $\widehat{\Gamma}$ of ${\mathbb Z}$ with respect to this group chain  admits a product decomposition into its Sylow $p$-subgroups
\begin{equation}\label{eq-pqlimit}
\widehat{\Gamma} ~ \cong ~ \prod_{i =1}^{\infty} \  {\mathbb Z}/q_i^{n(q_i)} {\mathbb Z} ~   \cdot ~ \prod_{p \in \pi_{\infty}} ~  \widehat{\mathbb Z}_{(p)} \ ,
\end{equation}
where $\widehat{\mathbb Z}_{(p)}$ denotes the $p$-adic completion of ${\mathbb Z}$.
Thus   $\pi(\xi(\widehat{\Gamma})) = \pi_f \cup \pi_{\infty}$. As ${\mathbb Z}$ is abelian,    $X_{\infty} = \widehat{\Gamma}$.
The type of this action classifies it up to return equivalence.
}
\end{ex}

 \begin{ex}\label{ex-almosttoral}
 {\rm 
The diagonal ${\mathbb Z}^n$-odometer actions   for $n \geq 2$ are direct extensions of Example~\ref{ex-toral}. 
Make $n$ choices of prime spectra as in Example~\ref{ex-toral}, then take the product action on the individual factors. The type of the action no longer determines the isomorphism class of the ${\mathbb Z}^n$-odometer actions obtained. The typeset is an invariant under return equivalence by Corollary~\ref{cor-abelian}, but only in special cases does the typeset determine the isomorphism class of the action. The interested reader can consult the works \cite{Arnold1982a,AV1992,Butler1965,Fuchs2015,Mutzbauer1983}.
 }
 \end{ex}

\begin{ex} \label{ex-twistedabelian} 
{\rm 
We next give an example that illustrates that the  commensurable relation on typesets is optimal. We construct the simplest example which shows this, and it is clear that many more similar constructions are possible.

Let $\Gamma = {\mathbb Z}^2 \rtimes {\mathbb Z}_2$  be the semi-direct product of ${\mathbb Z}^2$ with the order $2$ group ${\mathbb Z}_2 = {\mathbb Z}/2{\mathbb Z}$, where the generator $\sigma \in {\mathbb Z}_2$ acts on ${\mathbb Z}^2$ by permuting the summands. Let $\Gamma' = {\mathbb Z}^2$ which is abelian.

Choose distinct primes $p, q > 1$. Define the subgroup chain in $\Gamma$ and $\Gamma'$ as follows
\begin{equation}
\Gamma_{\ell} =  \{(p^{\ell} k, q^{\ell} m, id) \mid (h,m) \in {\mathbb Z}^2\} \quad , \quad \Gamma'_{\ell} = \{(p^{\ell} k, q^{\ell} m) \mid (k,m) \in {\mathbb Z}^2\} \ .
\end{equation}

Note  that the resulting actions $(X_{\infty}, \Gamma, \Phi)$ and $(X'_{\infty}, \Gamma', \Phi')$ are return equivalent. 

  Observe that  $C'_{\ell} = \Gamma'_{\ell}$. On the other hand,  $\Gamma_{\ell}$ is not normal in $\Gamma$, and   we have 
 \begin{equation}
C_{\ell} = \{((pq)^{\ell} k, (pq)^{\ell} m, id) \mid (h,m) \in {\mathbb Z}^2\} \subset \Gamma_{\ell} \ .
\end{equation}
  
The actions of $\Gamma$ and $\Gamma'$ have the same type, with characteristic functions $\chi(p) = \chi(q) = \infty$, and all other values are zero. However, we have the typesets 
\begin{equation}
\Xi[X_{\infty}, \Gamma, \Phi] =  \{ [(pq)^{\infty}] \}  \quad {\rm and}  \quad \Xi[X'_{\infty}, \Gamma', \Phi'] =  \{ [p^{\infty}] , [q^{\infty}], [(pq)^{\infty}]\} \ ,
\end{equation}
  so the typeset is not invariant under return equivalence.

This example easily generalizes, where we take $\Gamma' = {\mathbb Z} \oplus \cdots \oplus {\mathbb Z}$ to be the direct sum of $n$ copies of ${\mathbb Z}$, and replace ${\mathbb Z}_2$ with any non-trivial  subgroup   $\Delta \subset {\rm Perm}(n)$ of the permutation group on $n$ elements. Let $\Gamma = \Gamma' \rtimes \Delta$ be the semi-direct product of $\Gamma'$ with $\Delta$. Choose the subgroup chain $\{\Gamma'_{\ell'}\}$ in  $\Gamma'$ as in Example~\ref{ex-almosttoral} and use the   chain $\{\Gamma_{\ell} = \Gamma'_{\ell} \times id \}$ in $\Gamma$. Then the resulting actions are return equivalent, and one can obtain a wide variety of finite typesets $\Xi[X'_{\infty}, \Gamma', \Phi']$ for the abelian action. If the action of $\Delta$ is transitive, then  we have that $\Xi[X_{\infty}, \Gamma, \Phi]$ consists of a single element. When $\Delta$ does not act transitively,   there is even more variation on the typesets of the two actions.
}
\end{ex}

\subsection{Nilpotent odometers}
 
 The next examples use the actions associated to a renormalization of a finitely generated group $\Gamma$.
 Many   finitely generated nilpotent groups admit a renormalization \cite{Cornulier2016,DL2003a,DD2016,ER2005,LL2002,NekkyPete2011}, which yield many examples of odometer actions with well-defined typesets by Theorem~\ref{thm-renorm}.

 The integer Heisenberg group  is simplest non-abelian nilpotent group, and is  represented as the upper triangular matrices in ${\rm GL}(3,{\mathbb Z})$. That is,
\begin{equation}\label{eq-cH}
\Gamma =   \left\{  \left[ {\begin{array}{ccc}
   1 & a & c\\
   0 & 1 & b\\
  0 & 0 & 1\\
  \end{array} } \right] \mid a,b,c  \in {\mathbb Z}\right\} .
\end{equation}
 We denote a $3 \times 3$ matrix in $\Gamma$ by the coordinates as $(a,b,c)$.

 \begin{ex}\label{ex-trivial}
 {\rm
For a  prime $p \geq 2$, define the self-embedding $\varphi_p \colon \Gamma \to \Gamma$ by  
$\varphi(a,b,c) = (pa, pb, p^2c)$. Then define a group chain in $\Gamma$ by setting 
$$\Gamma_{\ell} = \varphi_p^{\ell}(\Gamma) = \{(p^{\ell} a, p^{\ell}b, p^{2\ell}c) \mid a,b,c \in {\mathbb Z}\} \quad, \quad \bigcap_{\ell > 0} \ \Gamma_{\ell} = \{e\} \ .$$

For $\ell > 0$, the normal core for $\Gamma_{\ell}$ is given by
$C_{\ell} = {\rm core}(\Gamma_{\ell})  = \{(p^{2\ell} a, p^{2\ell} b, p^{2\ell} c) \mid a,b,c \in {\mathbb Z}\}$, and so the quotient group  
$Q_{\ell} = \Gamma/C_{\ell} \cong \{( \overline{a}, \overline{b}, \overline{c}) \mid \overline{a}, \overline{b}, \overline{c} \in {\mathbb Z}/p^{2\ell}{\mathbb Z} \}$.
The profinite group $\widehat{\Gamma}_{\infty}$ is the inverse limit of the quotient groups $Q_{\ell}$ so we have
$\displaystyle \widehat{\Gamma}_{\infty} =   \{(\widehat{a}, \widehat{b}, \widehat{c}) \mid \widehat{a}, \widehat{b}, \widehat{c}\in \widehat{{\mathbb Z}}_{p^2} \}$.
   Thus,   every non-trivial $\gamma \in \Gamma$ has type $\tau[\gamma] = \tau[p^{\infty}]$.
 }
 \end{ex}
 
  \begin{ex}\label{ex-2primes}
 {\rm
For distinct  primes $p, q \geq 2$, define the self-embedding $\varphi_{p,q} \colon \Gamma \to \Gamma$ by  
$\varphi(a,b,c) = (pa, qb, pqc)$. Then define a group chain in $\Gamma$ by setting 
$$\Gamma_{\ell} = \varphi_{p,q}^{\ell}(\Gamma) = \{(p^{\ell} a, q^{\ell}b, (pq)^{\ell}c) \mid a,b,c \in {\mathbb Z}\} \quad, \quad \bigcap_{\ell > 0} \ \Gamma_{\ell} = \{e\} \ .$$

For $\ell > 0$, the normal core for $\Gamma_{\ell}$ is given by
$C_{\ell} = {\rm core}(\Gamma_{\ell})  = \{((pq)^{\ell} a, (pq)^{\ell} b, (pq)^{\ell} c) \mid a,b,c \in {\mathbb Z}\}$, and so we   obtain the quotient group  
$Q_{\ell} = \Gamma/C_{\ell} \cong \{( \overline{a}, \overline{b}, \overline{c}) \mid \overline{a}, \overline{b}, \overline{c} \in {\mathbb Z}/(pq)^{\ell}{\mathbb Z} \}$.
The profinite group $\widehat{\Gamma}_{\infty}$ is the inverse limit of the quotient groups $Q_{\ell}$ so we have
$\displaystyle \widehat{\Gamma}_{\infty} =   \{(\widehat{a}, \widehat{b}, \widehat{c}) \mid \widehat{a}, \widehat{b}, \widehat{c}\in \widehat{{\mathbb Z}}_{pq} \}$.
   Thus,   every non-trivial $\gamma \in \widehat{\Gamma}_{\infty}$ has type $\tau[\gamma] = \tau[(pq)^{\infty}]$.

Note that the typeset for the odometer   defined by the $\varphi_{p,q}$-renormalization   equals the typeset for the abelian action   in Example~\ref{ex-twistedabelian}, but the two actions are clearly not return equivalent.
 }
 \end{ex}
 
  A  second source of examples for odometer actions of nilpotent groups uses the decomposition of a profinite nilpotent group into its prime localizations, a technique that is especially adapted to realizing a given collection of primes as the spectrum of such an action. Moreover, these actions can be constructed to have    special dynamical properties, as in \cite{HL2022}.   The construction is necessarily more complex than for   renormalizable actions, as there must be an infinite sequence of choices to make.  These ideas were developed in \cite[Section~9]{HL2019} for actions of $SL(n, {\mathbb Z})$, and 
 the   work \cite{HL2023} discusses this construction for nilpotent groups.

  \section{Examples: $d$-regular actions} \label{sec-regular}

It is well-known that every odometer $({\mathfrak{X}},\Gamma,\Phi)$ is isomorphic to an action of $\Gamma$ on the boundary of a rooted tree. The study of actions on trees, especially the actions on $d$-ary (or $d$-regular) trees satisfying an additional condition of \emph{self-similarity}, is an active topic in Geometric Group Theory, see \cite{Nekrashevych2005,Grigorchuk2011} for surveys. The tree models for odometers are especially useful for constructing actions which are not topologically free, and thus for illustrating the commensurable relationship between types.  In this section, we study the typesets for odometers on boundaries of $d$-regular rooted trees.

\subsection{Actions on trees}\label{subsec-trees}

A \emph{tree} $T$ is an infinite graph with the set of vertices $V = \bigsqcup_{\ell \geq 0} V_\ell$ and the set of edges $E$. Each $V_\ell$, $\ell \geq 0$ is a finite set, called the \emph{set of vertices at level} $\ell$. Edges in $E$ join pairs of vertices in consecutive level sets $V_{\ell+1}$ and $V_\ell$, $\ell \geq 0$, so that a vertex in $V_{\ell+1}$ is connected to a single vertex in $V_\ell$ by a single edge. A tree is \emph{rooted} if $|V_0| = 1$.  

\begin{defn}
A tree $T$ is \emph{spherically homogeneous} if there is a sequence $n=(n_1,n_2,\ldots)$, called the \emph{spherical index} of $T$, such that for every $\ell \geq 1$ a vertex in $V_{\ell-1}$ is joined by edges to precisely $n_\ell$ vertices in $V_\ell$. In addition, $T$ is \emph{$d$-ary}, or \emph{$d$-regular}, if its spherical index $n=(n_1,n_2,\ldots)$ is \emph{constant}, that is, $n_\ell = d$ for some positive integer $d$. 
\end{defn}

We assume that $n_\ell \geq 2$ for $\ell \geq 1$. If $d = 2$, then a $2$-ary tree $T$ is also called a \emph{binary} tree.

Let $({\mathfrak{X}},\Gamma,\Phi)$ be an odometer, and let ${\mathcal U} = \{U_{\ell} \subset {\mathfrak{X}}  \mid \ell > 0\}$ be a choice of  an adapted neighborhood basis. By Corollary \ref{cor-Uchain} there is a group chain 
$\displaystyle  {\mathcal G}_{{\mathcal U}} = \{\Gamma_{\ell} = \Gamma_{U_{\ell}} \mid \ell \geq 0\}$ such that, associated to ${\mathcal G}_{{\mathcal U}}$ is an odometer   $(X_{\infty},\Gamma,\Phi_{\infty})$ which is  isomorphic  to $({\mathfrak{X}},\Gamma,\Phi)$. Here $X_\infty$ is the inverse limit space of finite sets $X_\ell = \Gamma/\Gamma_\ell$, given by \eqref{eq-invlimspace}, see Section \ref{subsec-gchains} for details.

A tree model for the action $({\mathfrak{X}},\Gamma,\Phi)$ is constructed using the group chain ${\mathcal G}_{{\mathcal U}}$.
For $\ell \geq 0$, let $V_\ell = X_\ell$, and join $v_\ell \in V_\ell$ and $v_{\ell+1} \in V_{\ell+1}$ by an edge if and only if $v_{\ell+1} \subset v_\ell$ as cosets. The  tree so constructed is spherically homogeneous, with spherical index entries $n_\ell = |\Gamma_{\ell-1}: \Gamma_{\ell}|$, for $\ell \geq 1$. The boundary $\partial T$ of $T$ is the collection of all infinite paths in $T$, that is,
  $$\partial T = \{(v_\ell)_{\ell \geq 0} \subset \prod_{\ell \geq 0} V_\ell \mid v_{\ell+1} \textrm{ and }v_\ell \textrm{ are joined by an edge}\} \cong X_\infty,$$
and the induced action of $\Gamma$ on $\partial T$, which we also denote by $(X_\infty,\Gamma,\Phi)$.

\begin{defn}
An odometer    $({\mathfrak{X}},\Gamma,\Phi)$ is $d$-\emph{regular}, or just \emph{regular}, if there exists $d \geq 2$ such that $({\mathfrak{X}},\Gamma,\Phi)$ is isomorphic to an action of $\Gamma$ on a rooted $d$-ary tree.
\end{defn}

\begin{remark}
{\rm 
An odometer    $({\mathfrak{X}},\Gamma,\Phi)$ is $d$-regular if the group chain $\mathcal{G}_{{\mathcal U}}$ above can be chosen so that each subgroup index $|\Gamma_{\ell}:\Gamma_{\ell -1}| = d$, for some $d \geq 2$ and all $\ell \geq 1$. Nilpotent actions given by a self-embedding $\varphi_p : \Gamma \to \Gamma$ in Example \ref{ex-trivial} are $d$-regular with $d = p^4$, and those in Example \ref{ex-2primes} are $d$-regular with $d = p^2q^2$.
}
\end{remark}

\subsection{Typesets of regular actions} We consider next the typeset of a $d$-regular odometer.

\begin{thm}\label{thm-regulartree}
An $d$-regular odometer   $({\mathfrak{X}},\Gamma,\Phi)$  has finite typeset $\Xi[\mathfrak{X},\Gamma,\Phi]$. 
More precisely, suppose $({\mathfrak{X}},\Gamma,\Phi)$ is isomorphic to an action on a $d$-ary rooted tree, for some $d \geq 2$. Let $P_d$ be the set of distinct prime divisors of the elements in the collection $\{2,3,\ldots,d\}$, and let $N_d = |P_d|$. Then
  \begin{align} \label{eq-upperboundtset}
  |\Xi[{\mathfrak{X}},\Gamma,\Phi]| \ \leq  \  \sum_{k=0}^{N_d} \binom{N_d}{k} = \sum_{k=0}^{N_d} \frac{N_d!}{(N_d - k)! k!} \ .
  \end{align}
  Moreover, each equivalence class   $\tau \in \Xi[{\mathfrak{X}},\Gamma,\Phi]$ is represented by   a Steinitz number $\xi$ with empty finite prime spectrum, $\pi_f(\xi) = \emptyset$, and so $\pi(\xi) = \pi_{\infty}(\xi)$. 
\end{thm}

\proof
Let $N_d$ and $P_d$ be as in the statement of the theorem, and let   $L_d = {\rm lcm} \{P_d\}$.

Let     ${\mathcal U} = \{U_{\ell} \subset {\mathfrak{X}}  \mid \ell > 0\}$ be a choice of  an adapted neighborhood basis with associated group chain 
$\displaystyle  {\mathcal G}_{{\mathcal U}} = \{\Gamma_{\ell} = \Gamma_{U_{\ell}} \mid \ell \geq 0\}$  and  $\Gamma$-odometer   $(X_{\infty},\Gamma,\Phi_{\infty})$  isomorphic to an action of $\Gamma$ on the boundary $\partial T$ of a  rooted  tree $T$ as constructed above.

\begin{lemma}\label{lemma-pdivides}
  Let $ \gamma \in \Gamma$, let $\xi(\gamma)$ be the Steinitz order of $\gamma$ as defined in Definition \ref{def-typegamma}, and let $p$ be such that $\chi_\xi(p) \ne 0$. Then $p$ divides $L_d$.
\end{lemma}

\proof
If $\chi_\xi(p) \ne 0$, then there exists the smallest $\ell \geq 1$ such that $p$ divides the order of the group $\langle \gamma \rangle /\langle \gamma \rangle_\ell$. The group $\langle \gamma \rangle /\langle \gamma \rangle_\ell  = \langle \gamma \rangle /\langle \gamma \rangle \cap C_\ell$ is isomorphic to a subgroup of $\Gamma/C_\ell$, where $C_\ell$ is the normal core of $\Gamma_\ell$, and acts on the coset space $X_\ell$ by permutations. Let $\lambda_{\gamma,\ell}$ be the permutation of $X_\ell$ induced by $\gamma$. Then the order of $\langle \gamma \rangle /\langle \gamma \rangle_\ell$ is equal to the order of $\lambda_{\gamma,\ell}$, and so equal to the least common multiple of the length of the cycles in $\lambda_{\gamma,\ell}$. 

Similarly, the order of $\langle \gamma \rangle /\langle \gamma \rangle_{\ell-1}$ is equal to the least common multiple of the length of the cycles in the permutation $\lambda_{\gamma,\ell-1}$ of $X_{\ell-1}$ induced by the action of $\gamma$. By the choice of $\ell$ the order of $\langle \gamma \rangle /\langle \gamma \rangle_{\ell-1}$ is not divisible by $p$. Therefore, for any cycle $c_{\ell-1}$ in $\lambda_{\gamma,\ell-1}$, the length $|c_{\ell-1}|$ is not divisible by $p$.

Consider the preimage $S_{c_{\ell-1}}$ of the set of elements in $c_{\ell-1}$ under the inclusion of cosets $X_\ell \to X_{\ell-1}$. Then $|S_{c_{\ell-1}}| = d|c_{\ell-1}|$, and $\gamma$ permutes the elements in $S_{c_{\ell-1}}$. Let $c_\ell$ be a cycle in the permutation $\mu_{c_{\ell-1}}$ of $S_{c_{\ell-1}}$ induced by $\gamma$. Since the action of $\gamma$ commutes with coset inclusions, $|c_\ell| = \alpha |c_{\ell-1}|$ for some $1 \leq \alpha \leq d$. Then $p$ must divide such an $\alpha$ for one of the cycles in $\mu_{c_{\ell-1}}$, for some cycle $c_{\ell-1}$ in $\lambda_{\gamma, \ell-1}$. It follows that $p$ divides $L_d$.
\endproof

Let $\xi(\gamma)$ be the Steinitz order of $\gamma$. Since every $p$ for which $0 < \chi_\xi(p) < \infty$, divides $L_d$ by Lemma \ref{lemma-pdivides}, then the finite prime spectrum $\pi(\xi(\gamma))$ is finite, and the type $\tau(\gamma)$ has a representative $\widehat{\xi}$, such that if $\chi_{\widehat{\xi}}(p) \ne 0$ then $\chi_{\widehat{\xi}}(p) = \infty$. It follows that two types $\tau(\gamma)$ and $\tau(\gamma')$ with respective Steinitz orders $\xi$ and $\xi'$ are distinct if and only if there exists a prime $p$ such that $\chi_\xi(p) = 0$ and $\chi_{\xi'}(p) \ne 0$. The bound in \eqref{eq-upperboundtset}, which is the number of distinct collections of prime divisors of $L_d$, follows.
\endproof

\begin{ex}
{\rm
Let $d = 2$, then $L_2 = 2$ and $N_2 = 1$, and the upper bound on the cardinality of the typeset for $2$-regular odometers is $2$, with possible types $\{[1], [2^\infty]\}$, where $1$ denotes the type of the identity element.

Let $d = 3$, then $L_3 = 6$ and $N_3 = 2$. Then the upper bound on the cardinality of the typeset for $3$-regular odometers is $4$, with possible types $\{[1], [2^\infty],[3^\infty],[(2 *3)^\infty]\}$.

Let $d = 4$, then $P_4 = \{2,3\}$, $L_4 = 6$ and $N_4 = 2$. Then the upper bound on the cardinality of the typeset for $4$-regular odometers is $4$, with possible types $\{[1], [2^\infty],[3^\infty],[(2 *3)^\infty]\}$.
}
\end{ex}

 \begin{ex}\label{ex-GuptaSidki}
  {\rm
  Let $p$ be an odd prime. The Gupta-Sidki $p$-group $GS(p)$ \cite[Section 1.8.1]{Nekrashevych2005}  is a group acting on the rooted $p$-ary tree. The group is generated by a cyclic permutation of the set of $p$ elements $\sigma = (0,1,\ldots,p-1)$, and the recursively defined map
    $$\tau = (\sigma, \sigma^{-1},1, \ldots, \tau).$$
 Every element in the Gupta-Sidki group $p$-group has finite order, and so $\Xi[X_\infty,GS(p),\Phi_{\infty}] = \{[1]\}$, where $[1]$ is the trivial type, i.e. the type of the identity element.
 }
  \end{ex}

 \begin{ex}\label{ex-alltypes}
{\rm
Let $d = 2$, and let $({\mathfrak X},\Gamma,\Phi)$ be a $2$-regular odometer action. Then $|\Xi[{\mathfrak X},\Gamma,\Phi]| \leq 2$, and either $\Xi[{\mathfrak X},\Gamma,\Phi] = \{[1]\}$, or $\Xi[{\mathfrak X},\Gamma,\Phi] = \{[1],[2^\infty]\}$. 

An example of an action where the typeset $\Xi[{\mathfrak X},\Gamma,\Phi] = \{[1]\}$, is the action of the Grigorchuk group,  see for instance \cite{Grigorchuk2011}. This group is an example of a Burnside group, which is an infinite group where every element has finite order, and so it has trivial typeset.

For the case $\Xi[{\mathfrak X},\Gamma,\Phi] = \{[1],[2^\infty]\}$ we have two classes of examples. First, the examples of odometers where $\tau(\gamma) = [2^\infty]$ for any non-trivial $\gamma \in \Gamma$, are the actions of iterated monodromy groups of quadratic polynomials for which the orbit of the critical point is periodic. Such groups are torsion-free, see \cite{BN2008}, and so they do not have any finite order element except the identity in $\Gamma$. Such actions include, but are not limited to, the action of the odometer on the binary tree, and the Basilica group.

Second, the action may have non-trivial element which have trivial type, and also non-trivial elements with type $\{[2^\infty]\}$. These are given by the actions of the iterated monodromy groups of quadratic polynomials with strictly pre-periodic orbits of the critical point, see \cite{BN2008} for details.}
\end{ex}

\subsection{Typeset under the commensuration relation}
   
 We now show that $H$-restricted types of a group element and of its conjugate need not coincide.
 
 Given a tree $T$ and a vertex $v \in V$, denote by $T_v$ the subtree of $T$ with root $v$.
 
 \begin{defn}\label{defn-weaklybranch}\cite{Grigorchuk2011}
 Let $T$ be a spherically homogeneous tree, and let $\Gamma \subset {\rm Aut}(T)$.  The action of $\Gamma$ on $\partial T$ is \emph{weakly branch} if the following conditions hold:
 \begin{enumerate}
 \item The restriction of the action of $\Gamma$ to each vertex level set $V_n$, $n \geq 1$, is transitive.
 \item For every $\ell \geq 1$ and every $v \in V_\ell$, there exists $g \in \Gamma$ such that $g \cdot v = v$, the restriction of the action of $g$ on the subtree $T_v$ is non-trivial, and the restriction of $g$ to the complement $\partial T - \partial T_v$ is the identity map.
 \end{enumerate}
 \end{defn}
   
 \begin{ex}\label{ex-branch}
 {\rm 
Let $T$ be a spherically homogeneous tree and let $\Gamma \subset {\rm Aut} (T)$ be so that the action of $\Gamma$ on $\partial T$ is weakly branch.  Choose a sequence $(v_\ell)$ of vertices in $T$, and let $\Gamma_\ell$ be the isotropy group (equivalently, the stabilizer) of $v_\ell \in V_\ell$. Then $\Gamma_\ell$ fixes $v_\ell$ while permuting other vertices in $V_\ell$.

Choose $n \geq 1$. The normal core $C_n$ of $\Gamma_n$ consists of elements of $\Gamma_n$ which fix every vertex in $V_n$. Let $g \in \Gamma_n$ be an element given by Definition \ref{defn-weaklybranch} of a weakly branch group, namely, the restriction of $\gamma$ to $T_{v_n}$ is non-trivial, and $\gamma$ is trivial on the complement of $\partial T_{v_n}$ in $\partial T$. In particular, this means that $\gamma \in C_n$. However, since $\gamma$ acts non-trivially on $\partial T_{v_n}$, then there exists $\ell > n$ such that $\gamma \notin C_\ell$, and then the Steinitz order $\xi(\gamma)$ is non-trivial.

Now let $\delta \in \Gamma$ be such that $\delta \cdot v_n \ne v_n$. Then the restriction of $\delta \gamma \delta^{-1}$ to $\partial T_{v_n}$ is the identity map. Since the Steinitz order is invariant under conjugation, $\xi (\delta \gamma \delta^{-1})$ is non-trivial.

Now consider the $\Gamma_n$-restricted Steinitz orders of $\gamma$ and of $\delta \gamma \delta^{-1}$. The subgroup $C^{\Gamma_n}_{\ell}$ is the normal core of $\Gamma_\ell$ in $\Gamma_n$, so it fixes every vertex in the set $V_\ell \cap T_{v_n}$, and it may permute the vertices of $V_\ell$ which are not in $T_{v_n}$. In particular, for all $\ell \geq n$ we have $\delta \gamma \delta^{-1} \in C^{\Gamma_n}_\ell$, and the $\Gamma_n$-restricted Steinitz order $\xi^{\Gamma_n} (\delta \gamma \delta^{-1})$ is trivial. At the same time,  $\xi^{\Gamma_n}(\gamma)$ is non-trivial, since the action of $\gamma$ permutes the vertices in $V_\ell \cap T_{v_n}$ for some $\ell \geq n$. Depending on the group $\Gamma$, $\xi^{\Gamma_n}(\gamma)$ may have the trivial or non-trivial type; it is straightforward to construct examples of both situations.
 }
 \end{ex}

 \section{Solenoidal manifolds}\label{sec-solenoids}
 
Solenoidal manifolds were introduced in Section~\ref{subsec-solenoids}. In this section, we discuss the relation between the classification of solenoidal manifolds up to homeomorphism, and the classification of odometers up to return equivalence. We then give the proof of Theorems~\ref{thm-main11} and \ref{thm-main21}.
 
Recall that a \emph{presentation} is a  sequence of \emph{proper finite covering} maps 
${\mathcal P} = \{\, q_{\ell} \colon  M_{\ell} \to M_{\ell -1} \mid  \ell \geq 1\}$, where each $M_{\ell}$ is a compact connected manifold without boundary of dimension $n$.   The inverse limit     
\begin{equation}\label{eq-presentationinvlim2}
{\mathcal M}_{{\mathcal P}} \equiv \lim_{\longleftarrow} ~ \{ q_{\ell } \colon M_{\ell } \to M_{\ell -1}\} ~ \subset \prod_{\ell \geq 0} ~ M_{\ell} ~  
\end{equation}
is   the  \emph{weak solenoid}, or \emph{solenoidal manifold},    associated to ${\mathcal P}$.  

 By the definition of the inverse limit, for a sequence $\{x_{\ell} \in M_{\ell} \mid \ell \geq 0\}$, we have 
\begin{equation}\label{eq-presentationinvlim3}
x = (x_0, x_1, \ldots ) \in {\mathcal M}_{{\mathcal P}}   ~ \Longleftrightarrow  ~ q_{\ell}(x_{\ell}) =  x_{\ell-1} ~ {\rm for ~ all} ~ \ell \geq 1 ~. 
\end{equation}
For each $\ell \geq 0$, there  is a    fibration ${\widehat{q}}_{\ell} \colon {\mathcal M}_{{\mathcal P}} \to M_{\ell}$,  given by projection onto the $\ell$-th factor in \eqref{eq-presentationinvlim2}. Also,  there is a covering map  denoted by
$\overline{q}_{\ell} = q_{\ell} \circ q_{\ell -1} \circ \cdots \circ q_1 \colon M_{\ell} \to M_0$, such  that ${\widehat{q}}_0 = \overline{q}_{\ell} \circ {\widehat{q}}_{\ell}$. 
Choose a basepoint $x_0 \in M_0$ and basepoint $x_{\infty} \in {\mathfrak{X}} = {\widehat{q}}_0^{-1}(x_0)$, the fiber over $x_0$. Then for each $\ell > 0$, this defines the basepoint $x_{\ell} =  {\widehat{q}}_{\ell}(x_{\infty}) \in M_{\ell}$.
Let $\Gamma_{\ell} = (\overline{q}_{\ell})_{\#}(\Gamma_0)$ denote the image $\Gamma_{\ell}$ of the fundamental group $\pi_1(M_{\ell}, x_{\ell})$ in $\Gamma = \pi_1(M_0, x_0)$.

Let ${\mathcal G}_{\mathcal P} = \{\Gamma_{\ell} \mid \ell \geq 0\}$ be the group chain determined by $\mathcal P$ and the choice of basepoint $x \in {\mathcal M}_{{\mathcal P}}$.
Let $(X_{\infty},\Gamma,\Phi_{\infty})$ denote the odometer defined in \eqref{eq-invlimspace} which is determined by ${\mathcal G}_{\mathcal P}$ which  is called the monodromy action of the fibration ${\widehat{q}}_0 \colon {\mathcal G}_{\mathcal P} \to M_0$. Note that we suppress the dependence on the choice of the basepoint $x$ for convenience of notation.

The  fiber $\overline{q}_{\ell}^{-1}(x_0) \subset M_{\ell}$ is 
 identified with the quotient set $X_{\ell} = \Gamma/\Gamma_{\ell}$ which is a left $\Gamma$-space. In this way, the subspace ${\mathfrak{X}} = {\widehat{q}}_0^{-1}(x_0)$ of the inverse limit in \eqref{eq-presentationinvlim2} is identified with the $\Gamma$-space $X_{\infty}$. We then have the following result:
   \begin{thm}  \cite[Theorem~1.1]{CHL2019} \label{thm-morita}
  Let  ${\mathcal M}_{{\mathcal P}}$ and  ${\mathcal M}_{{\mathcal P}'}$ be homeomorphic solenoidal manifolds. Then the corresponding  odometers $(X_{\infty},\Gamma,\Phi_{\infty})$ and $(X_{\infty}',\Gamma',\Phi_{\infty}')$ are return equivalent. 
\end{thm}
Thus,  an invariant of odometers modulo return equivalence  is an invariant for   solenoidal manifolds modulo homeomorphism.  There are various classes of solenoidal manifolds where the converse statement can be shown, as   discussed in \cite{CHL2019}. Here is one example:
\begin{thm} \cite[Theorem~1.3]{CHL2019}\label{thm-toroidal}
Let  ${\mathcal M}_{{\mathcal P}}$ and  ${\mathcal M}_{{\mathcal P}'}$ be  toroidal  solenoidal of the same dimension; that is, both presentations ${\mathcal P}$ and ${\mathcal P}'$ have  base manifold the $n$-torus ${\mathbb T}^n$. If the monodromy actions $(X_{\infty},\Gamma,\Phi_{\infty})$ and $(X_{\infty}',\Gamma',\Phi_{\infty}')$ are return equivalent, then  ${\mathcal M}_{{\mathcal P}}$ and  ${\mathcal M}_{{\mathcal P}'}$ are homeomorphic.
\end{thm}

 \subsection{Type invariants for solenoidal manifolds}\label{subsec-typesolenoids}
 
 The covering degree $m_{\ell}$ of $q_{\ell} \colon  M_{\ell} \to M_{\ell -1}$ equals the index of the subgroup $[\Gamma_{\ell-1} : \Gamma_{\ell}]$, and so   the covering degree of $\overline{q}_{\ell}   \colon M_{\ell} \to M_0$ is given by 
\begin{equation}
{\rm deg}(\overline{q}_{\ell}) = m_{\ell} \cdot m_{\ell -1} \cdots m_1 = [\Gamma : \Gamma_{\ell}].
\end{equation}
The Steinitz order of a presentation $\mathcal P$ is 
   \begin{equation}\label{eq-highersteinitzorder33}
\xi({\mathcal P}) = {\rm lcm} \{ m_1   m_2 \cdots m_{\ell} \mid  \ell > 0\} = {\rm lcm} \{ [\Gamma : \Gamma_{\ell}] \mid  \ell > 0\}  = \xi({\mathcal G}_{\mathcal P}) \ .
\end{equation}
 The Steinitz number $\xi({\mathcal P}')$ can be thought of as the covering degree of the fibration ${\widehat{q}}_0 \colon {\mathcal M}_{{\mathcal P}} \to M_0$.

Given a second presentation ${\mathcal P}'$ with solenoidal manifold ${\mathcal M}'_{{\mathcal P}'}$ and choice of basepoint $x_0' \in {\mathcal M}'_{{\mathcal P}'}$, we similarly obtain subgroups $\Gamma'_{\ell}$ defining the group chain ${\mathcal G}'_{{\mathcal U}'}$, and we have $\xi({\mathcal P}')   = \xi({\mathcal G}'_{{\mathcal U}'})$.

Suppose that the solenoidal manifolds ${\mathcal M}_{{\mathcal P}}$ and ${\mathcal M}'_{{\mathcal P}'}$ are homeomorphic, then by Theorem~\ref{thm-morita}, the corresponding  odometers $(X_{\infty},\Gamma,\Phi_{\infty})$ and $(X_{\infty}',\Gamma',\Phi_{\infty}')$ are return equivalent, and hence they have the equal types, $\tau[X_{\infty},\Gamma,\Phi_{\infty}] = \tau[X_{\infty}',\Gamma',\Phi_{\infty}']$  by Theorem~\ref{thm-main2}. We thus obtain:

  \begin{thm}\label{thm-main1111}  
 Associated to a presentation $\mathcal P$ is a well-defined Steinitz number $\xi({\mathcal P})$, and the  type $\tau[\xi({\mathcal P})]$    depends only on the homeomorphism class of ${\mathcal M}_{{\mathcal P}}$, which is denoted  by $\tau[{\mathcal M}_{{\mathcal P}}]$.
 \end{thm}

We thus obtain the following well-defined homeomorphism invariants for solenoidal manifolds:

\begin{cor}\label{cor-spectrum}  The \emph{infinite prime spectrum} $\pi_{\infty}(\xi({\mathcal P}))$ of a   solenoidal manifold  ${\mathcal M}_{{\mathcal P}}$ 
   depends only on the homeomorphism class of ${\mathcal M}_{{\mathcal P}}$, as does also the property that $\pi_f(\xi({\mathcal P}))$ is  infinite.
 \end{cor}

For orientable 1-dimensional solenoids, where each map $q_{\ell} \colon {\mathbb S}^1 \to {\mathbb S}^1$ is orientable,   Bing observed in \cite{Bing1960} that    if  $\tau[{\mathcal M}_{{\mathcal P}}] = \tau[{\mathcal M}_{{\mathcal P}'}]$ then   ${\mathcal M}_{\mathcal P}$ and ${\mathcal M}_{{\mathcal P}'}$  are homeomorphic.  McCord showed in \cite[Section~2]{McCord1965}   the converse, that if  ${\mathcal M}_{\mathcal P}$ and ${\mathcal M}_{{\mathcal P}'}$  are homeomorphic, then  $\tau[{\mathcal M}_{{\mathcal P}}] = \tau[{\mathcal M}_{{\mathcal P}'}]$.   (See also   \cite{AartsFokkink1991}.)  Together these results yield the well-known classification:  
\begin{thm}\label{thm-onedimSol2}
For orientable 1-dimensional solenoidal manifolds, 
  ${\mathcal M}_{\mathcal P}$ and ${\mathcal M}_{{\mathcal P}'}$ are homeomorphic  if and only if    $\tau[{\mathcal M}_{\mathcal P}] = \tau[{\mathcal M}_{{\mathcal P}'}]$. 
\end{thm}
For solenoidal manifolds of   dimension $n \geq 2$, no such classification by invariants can exist, even when the base manifold $M_0 = {\mathbb T}^n$. The most one can hope for is to define invariants which distinguish various homeomorphism classes, such as the type and typeset for a solenoidal manifold, invariants  derived from the Cantor action associated to the monodromy of the fibration map ${\widehat{q}}_0 \colon {\mathcal M}_{{\mathcal P}} \to M_0$.  

Note that by Theorem~\ref{thm-toroidal}, if two toroidal solenoids have return equivalent monodromy odometers, then they are homeomorphic. Thus, in this special case, the classification problem for solenoidal manifolds is equivalent to the classification problem of ${\mathbb Z}^n$-odometers modulo return equivalence. Giordano, Putnam and Skau discuss in \cite{GPS2019} this classification problem.

 \subsection{Typeset invariants for solenoidal manifolds}\label{subsec-typesetsolenoids}

The typeset invariants for solenoidal manifolds provide more  refined invariants of their homeomorphism class. Before stating the precise definition, we give an intuitive definition.

Let $M_0$ be the base manifold of a presentation ${\mathcal P}$.
Let $x_{\infty} \in {\mathcal M}_{{\mathcal P}}$ be a choice of a basepoint, which determines the basepoint  $x_0  \in M_0$ and fundamental group  $\Gamma = \pi_1(M_0, x_0)$. Suppose that $\gamma \in \Gamma$     is represented by a simple closed curve  $c_{\gamma} \colon {\mathbb S}^1 \to M_0$. 
The projection map ${\widehat{q}}_0 \colon {\mathcal M}_{{\mathcal P}} \to M_0$   restricts to a covering map on each leaf of the foliation ${\mathcal F}_{{\mathcal P}}$ of ${\mathcal M}_{{\mathcal P}}$. When the preimage ${\mathcal M}_{\gamma} = {\widehat{q}}_0^{-1}(c_{\gamma}({\mathbb S}^1))$ is connected, it is a 1-dimensional solenoid, and so has a well-defined type by Theorem~\ref{thm-main1111}. Denote this type by   $\tau[\gamma]$. The collection of these types is the `` intuitive typeset'' of ${\mathcal M}_{{\mathcal P}}$.

 The formal definition of the typeset proceeds as in Section~\ref{subsec-typesolenoids}.  Associate to a presentation ${\mathcal P}$ the monodromy odometer $(X_{\infty}, \Gamma, \Phi_{\infty})$. Then for each $\gamma \in \Gamma$ there is the type 
  $\tau[\gamma]$ as in Definition~\ref{def-typegamma}.

\begin{defn}\label{def-typesetP}  The   \emph{typeset} of  ${\mathcal P}$ is the countable collection of types   $\Xi[{\mathcal P}] = \{  \tau[\gamma] \mid \gamma \in \Gamma\}$.
 \end{defn}

\begin{thm}\label{thm-main211}  Suppose that the   solenoidal manifolds ${\mathcal M}_{{\mathcal P}}$ and ${\mathcal M}_{{\mathcal P}'}$ are homeomorphic, then their typesets $\tau[{\mathcal P}]$ and $\tau[{\mathcal P}']$ are commensurable. 
 \end{thm}
\proof
Suppose that ${\mathcal M}_{{\mathcal P}}$ and ${\mathcal M}_{{\mathcal P}'}$ are homeomorphic, then by Theorem~\ref{thm-main1111} their monodromy odometers $(X_{\infty}, \Gamma, \Phi_{\infty})$ and $(X_{\infty}', \Gamma', \Phi'_{\infty})$ are return equivalent. Then by Theorem~\ref{thm-main2b} the typesets $\Xi[{\mathcal P}]$ and $\Xi[{\mathcal P}']$ are commensurable.
\endproof

 The type and typeset invariants suggest numerous questions about their relation to the analytic and geometric properties of solenoidal manifolds \cite{Sullivan2014,Verjovsky2014,Verjovsky2022}.

  
 \bibliographystyle{plain}

\end{document}